\documentclass[12pt]{amsart}
\usepackage{amssymb}
\usepackage[noadjust]{cite}
\usepackage{booktabs}
\usepackage{mathrsfs}
\usepackage{url}
\usepackage{hyphenat}
\usepackage{mathtools}
\usepackage{cancel} 
\usepackage[T1]{fontenc}
\usepackage[utf8]{inputenc}
\usepackage{mdwtab}
\usepackage[pdftex,lmargin=1in,rmargin=1in,tmargin=1in,bmargin=1in]{geometry}
\usepackage[bookmarks=true, bookmarksopen=true,%
    bookmarksdepth=3,bookmarksopenlevel=2,%
    colorlinks=true,%
    linkcolor=blue,%
    citecolor=blue,%
    filecolor=blue,%
    menucolor=blue,%
    urlcolor=blue]{hyperref}
\hypersetup{pdftitle={Faithful linear-categorical mapping class group action}}
\hypersetup{pdfauthor={Robert Lipshitz, Peter S. Ozsváth, and Dylan P. Thurston}}
\usepackage{color}
\usepackage{tikz}

\usetikzlibrary{matrix,arrows}
\tikzstyle{every picture}=[> = latex']
\tikzset{cdlabel/.style={above,sloped,
    execute at begin node=$\scriptstyle,execute at end node=$}}
\tikzset{algarrow/.style={->, thick}}
\tikzset{blgarrow/.style={->, thick}}
\tikzset{clgarrow/.style={->, thick}}
\tikzset{tensoralgarrow/.style={->, thin, double}}
\tikzset{tensorblgarrow/.style={->, thin, double}}
\tikzset{tensorclgarrow/.style={->, thin, double}}
\tikzset{modarrow/.style={->, dashed}}
\tikzset{othmodarrow/.style={->, thick}}
\tikzset{Amodar/.style={->, dashed}}
\tikzset{Dmodar/.style={->, dashed}}

\newcommand{\RR}{\mathbb R}
\newcommand{\CC}{\mathbb C}
\newcommand{\ZZ}{\mathbb Z}

\newcommand{\FF}{\mathbb F}
\newcommand{\NN}{\mathbb N}

\newcommand{\bD}{\mathbb{D}}

\newcommand{\co}{\colon}

\newcommand{\bdy}{\partial}

\newcommand{\lbracket}{[}
\newcommand{\rbracket}{]}

\newcommand{\Hyph}{\text{-}}

\DeclareMathOperator{\spin}{spin}
\newcommand{\SpinC}{\spin^c}

\DeclareMathOperator{\supp}{supp}

\theoremstyle{plain}
\newtheorem{theorem}{Theorem}
\numberwithin{equation}{section}

\newtheorem{proposition}[equation]{Proposition}
\newtheorem{lemma}[equation]{Lemma}
\newtheorem{corollary}[equation]{Corollary}

\newtheorem{definition}[equation]{Definition}

\theoremstyle{definition}

\theoremstyle{remark}
\newtheorem{example}[equation]{Example}
\newtheorem{remark}[equation]{Remark}

\newcommand{\HF}{\mathit{HF}}

\newcommand{\CF}{\mathit{CF}}

\newcommand{\x}{\mathbf x}
\newcommand{\y}{\mathbf y}

\newcommand\HH{\mathit{HH}}

\newcommand\Hochschild\HH

\newcommand{\MCG}{\mathit{MCG}}

\newcommand{\Ainf}{\mathcal A_\infty}
\newcommand{\Alg}{\mathcal{A}}

\newcommand\Blg{\mathcal{B}}
\newcommand\Clg{\mathcal{C}}

\newcommand{\Idem}{\mathcal{I}}
\newcommand{\alphas}{{\boldsymbol{\alpha}}}
\newcommand{\betas}{{\boldsymbol{\beta}}}

\newcommand{\cM}{\mathcal{M}}

\newcommand{\DD}{\textit{DD}}
\newcommand{\DA}{\textit{DA}}

\newcommand{\AAm}{\textit{AA}} 

\newcommand{\CFDD}{\mathit{CFDD}}

\newcommand{\CFDA}{\mathit{CFDA}}
\newcommand{\CFDAa}{\widehat{\CFDA}}
\newcommand{\CFAA}{\mathit{CFAA}}
\newcommand{\CFAAa}{\widehat{\CFAA}}

\newcommand{\CFDDa}{\widehat{\CFDD}}

\newcommand{\cZ}{\mathcal{Z}}
\newcommand{\PtdMatchCirc}{\cZ}
\newcommand{\PMC}{\PtdMatchCirc}

\newcommand{\CircPts}{{\mathbf{a}}}
\newcommand{\bbpt}{\mathbf{z}}

\newcommand\Id{\mathbb{I}}

\newcommand\DTP{\mathbin{\widetilde\otimes}}
\newcommand\DT{\boxtimes}
\newcommand\Gen{\mathfrak{S}}

\newcommand{\Field}{{\FF_2}}

\newcommand{\Heegaard}{\mathcal{H}}
\newcommand{\HD}{\Heegaard}

\newcommand{\st}{^\text{st}}

\newcommand{\ModCat}{\mathsf{Mod}}

\newcommand{\Cat}{\mathscr{C}}

\DeclareMathOperator{\Mor}{Mor}

\newcommand{\DerBounded}{\mathcal{D}^b}

\newcommand{\op}{\mathrm{op}}
\newcommand\PunctF{F^\circ}

\makeatletter
\newcommand\honestalg[3]{\bigl\lbracket
\begin{smallmatrix} #1\@ifempty{#3}{}{&#3} \\ #2 \end{smallmatrix}
\bigr\rbracket}

\makeatother
\newcommand{\lab}[1]{$\scriptstyle #1$}

\newcommand{\lsub}[2]{{}_{#1}#2}

\newread\testin

\graphicspath{{draws/}{mpdraws/}}
\makeatletter
\def\input@path{{}{draws/}}
\makeatother

\def\mathcenter#1{%
  \vcenter{\hbox{$#1$}}%
}

\makeatletter
\newcommand\mi@kern[1]{%
  \settowidth\@tempdima{$\mi@obj^{#1}$}
  \kern-\@tempdima
  #1
  \settowidth\@tempdima{$\mi@obj$}
  \kern\@tempdima
}

\newtoks\mi@toksp
\newtoks\mi@toksb
\DeclareRobustCommand{\manyindices}[5]{
  \def\mi@obj{#5}
  \mi@toksp\expandafter{\mi@kern{#2}}
  \mi@toksb\expandafter{\mi@kern{#1}}
  \@mathmeasure4\textstyle{#5_{#1}^{#2}}
  \@mathmeasure6\textstyle{#5_{#3}^{#4}}
  \dimen0-\wd6 \advance\dimen0\wd4
  \@mathmeasure8\textstyle{\hphantom{{}_{#1}^{#2}}#5^{\the\mi@toksp#4}_{\the\mi@toksb#3}}
  \hbox to \dimen0{}{\kern-\dimen0\box8}
}
\makeatother 

\usepackage{ifpdf}
\ifpdf
  \let\textalt\texorpdfstring
\else
  \newcommand{\textalt}[2]{#1}
\fi

\newcommand{\DDHalfId}{\mathit{DD}\bigl({\textstyle \frac{\Id}{2}}\bigr)}
\newcommand{\Matching}{M}
\newcommand{\Diagram}{\mathcal{D}}
\newcommand{\DM}{P}
\newcommand{\ADM}{Q}
\newcommand{\AM}{M}
\newcommand{\DAM}{N}
\newcommand{\ShortChords}{\mathsf{SC}}

\begin{document}
\title[Faithful linear-categorical mapping class group action]
{A faithful linear-categorical action of the mapping class group of a surface with boundary}

\author[Lipshitz]{Robert Lipshitz}
\thanks{RL was supported by NSF grant DMS-0905796 and a Sloan Research
  Fellowship.}
\address{Department of Mathematics, Columbia University\\
  New York, NY 10027}
\email{lipshitz@math.columbia.edu}

\author[Ozsv\'ath]{Peter~S.~Ozsv\'ath}
\thanks{PSO was supported by NSF grant DMS-0505811 and a Clay
  Senior Scholar Fellowship.}
\address {Department of Mathematics, MIT\\ Cambridge, MA 02139}
\email {petero@math.mit.edu}

\author[Thurston]{Dylan~P.~Thurston}
\thanks {DPT was supported by NSF
  grant DMS-1008049.}
\address{Department of Mathematics,
         Barnard College,
         Columbia University\\
         New York, NY 10027}
\email{dthurston@barnard.edu}


\begin{abstract}
  We show that the action of the mapping class group on bordered Floer
  homology in the second to extremal $\SpinC$-structure is faithful. 
  This paper is designed partly as an introduction to the subject, and
  much of it should be readable without a background in Floer
  homology.
\end{abstract} 

\maketitle 

\tableofcontents
\section{Introduction}

Two long-standing, and apparently unrelated, questions in
low-dimensional topology are whether the mapping class group of a
surface is linear
and whether the Jones polynomial detects the unknot. In 2010,
Kronheimer-Mrowka gave an affirmative answer to a categorified version
of the second question: they showed that Khovanov homology, a
categorification of the Jones polynomial, does detect the unknot~\cite{KronheimerMrowka11:detect}.
(Previously, Grigsby and Wehrli had shown that any nontrivially-colored
Khovanov homology detects the unknot~\cite{GrigsbyWehrli:detects}.) In
this paper, we give an affirmative answer to a categorified version of
the first question. That
is, while we do not know if the mapping class group of a surface (with
boundary) acts faithfully on
a finite-dimensional linear space, we are able to give an explicit faithful action
on a finitely-generated linear (in fact,
triangulated) category.\footnote{Because in this paper we do not
  discuss gradings, which are somewhat subtle,
 the categories will actually be ungraded analogues of triangulated
 categories.  See, e.g.,
 \cite[Section~\ref*{LOT2:sec:algebras-gradings}]{LOT2} for more on
 the gradings in bordered Floer theory.} The decategorification of
this action is the standard
action of the mapping class group on $H_1$; see
Theorem~\ref{thm:grothendiek} in Section~\ref{sec:finite-generation}.

In more detail, the structure is as follows. To a surface $F$ with
boundary and a marked point on each boundary component, we associate a
finite-dimensional algebra $\Blg(F)$ over $\Field=\ZZ/2$. (There is
some choice in the definition of $\Blg(F)$; see
Section~\ref{sec:algebras}.) To a mapping class $\phi\co F\to F$,
fixing the boundary, we associate a quasi-isomorphism class of finite-dimensional differential
$\Blg(F)$-bimodules $\CFDAa(\phi)$. These have the property that
\begin{equation}\label{eq:DADA-pairing}
\CFDAa(\psi\circ\phi)\simeq\CFDAa(\phi)\otimes_{\Blg(F)}\CFDAa(\psi).\footnotemark
\end{equation}
\footnotetext{More
    honestly, here $\otimes$ should be interpreted as the derived, or $\Ainf$, tensor
    product, though it is possible to work with models for $\CFDAa$
    for which this agrees with the ordinary tensor product. Later, we
    will use $\DTP$ for the $\Ainf$-tensor product to keep track of
    this distinction.}
Moreover, 
\begin{align}
\CFDAa(\Id)&\simeq \lsub{\Blg(F)}\Blg(F)_{\Blg(F)},\label{eq:id-is-id}
\end{align}
where $\lsub{\Blg(F)}\Blg(F)_{\Blg(F)}$ denotes the algebra
$\Blg(F)$ viewed as a bimodule over itself.

Let $\lsub{\Blg(F)}\ModCat$ denote the category of finitely-generated
left $\Blg(F)$-modules. For each mapping class $\phi$ we have a functor
$\Phi_\phi\co \lsub{\Blg(F)}\ModCat\to \lsub{\Blg(F)}\ModCat$ given by
$\Phi_\phi(\cdot)=\CFDAa(\phi)\otimes_{\Blg(F)}\cdot$. Equations~(\ref{eq:DADA-pairing})
and~(\ref{eq:id-is-id}) almost imply that this is an action; the main
defect is that Equation~(\ref{eq:DADA-pairing}) only gives homotopy
equivalences, not isomorphisms (or equalities). To rectify this, we
replace $\lsub{\Blg(F)}\ModCat$ with the associated derived
category $\DerBounded(\lsub{\Blg(F)}\ModCat)$ of finitely-generated
modules. (This is quite
concrete: since finite-dimensional modules over our algebras admit
finite-dimensional projective resolutions, $\DerBounded(\lsub{\Blg(F)}\ModCat)$ is just the homotopy
category of finitely-generated projective modules over
$\Blg(F)$.)
Equations~(\ref{eq:DADA-pairing}) and~(\ref{eq:id-is-id}) then imply
that tensoring with the modules $\CFDAa(\psi)$ gives an action of the
mapping class group on $\DerBounded(\lsub{\Blg(F)}\ModCat)$.
(There are some subtleties related to group actions on categories. See
for example \cite[Section 8]{LOT2} for a review of the relevant
definitions.)

The bimodules $\CFDAa(\phi)$ carry geometric information. In particular,
the rank of the homology of $\CFDAa(\phi)$ is given by a certain
intersection number. This turns out to be enough
to prove that
\begin{equation}
\CFDAa(\phi) \not\simeq \CFDAa(\Id) \qquad \text{if }\phi \not\sim \Id.\label{eq:phi-not-id}
\end{equation}
As a corollary, we have:
\begin{theorem}\label{thm:intro-faithful}
  The action of the mapping class group $\MCG_0(F)$ on
  $\DerBounded(\lsub{\Blg(F)}\ModCat)$ given by tensoring with the
  bimodules $\CFDAa(\phi)$ is faithful.
\end{theorem}

In fact, there are two different ways
we can do this construction combinatorially. One leads to somewhat
simpler algebras, but more complicated ($\Ainf$) bimodules; the other
leads to more complicated (differential) algebras but simpler
(differential) bimodules. Although these two actions are equivalent in a
certain sense---see Proposition~\ref{prop:actions-equivalent},
below---we will give both approaches.

Experts in bordered Floer theory are warned that throughout this paper
we are working in the second to extremal $\SpinC$-structure. In the
notation of~\cite{LOT1}, the algebras $\Blg(F)$ (respectively
$\Clg(F)$) in this paper are $\Alg(F,\allowbreak -g+1)$ (respectively
$\Alg(F,\allowbreak g-1)$), where $g$ is the genus of $F$, and the
bimodules $\CFDAa(\phi)$ are the corresponding summands of the
bimodules $\CFDAa(\phi)$ from~\cite{LOT2}.

This paper has two main goals. The first goal is to prove faithfulness
of the mapping class group action (Theorem~\ref{thm:intro-faithful}).
The proof of faithfulness itself is short, and the reader familiar
with the bordered Floer package may wish to skip directly to
Section~\ref{sec:faithfulness} (perhaps after perusing some of the
pictures earlier in the paper), where the proof is given. The second
goal is to give a
combinatorial description of this mapping class group action (in the
second to extremal $\SpinC$-structure).  This paper is partly intended
as an introduction to the subject. So, we include a complete
description of the relevant algebras and modules. The proof of
faithfulness is also elementary, and both the modules and the
faithfulness proof are closely related to familiar tools in mapping
class group theory.  We do not give self-contained proofs that the
bimodules associated to mapping classes are well-defined, or that
tensoring with them gives a well-defined action; these results draw
on~\cite{LOT2},
which uses the theory of pseudoholomorphic curves. Since the
first version of this paper was written, Kyler Siegel has given direct
combinatorial proofs of these facts; see~\cite{Siegel:mcgAction}.

In this paper, we treat mapping class groups of any surface with
non-empty boundary. The case of actions of braid groups on
triangulated categories (unlike the more general case) has received
substantial attention in the literature. See in
particular~\cite{KhS02:BraidGpAction}, and
also~\cite{KhTomas07:CobordismsCategories} and the references
contained therein. Another triangulated category on which the mapping
class group acts is the Fukaya category of a surface; a theorem of
Seidel \cite[Theorem
1]{Seidel02:FloerMappingClass}, together with a folk conjecture
relating the Hochschild homology of functors on the Fukaya category
and Floer homology of symplectomorphisms,
should imply this action is faithful for a closed surface. 
The argument in Section~\ref{sec:faithfulness}, which was
inspired by \cite{KhS02:BraidGpAction}, can be adapted to give a more
direct proof of faithfulness of the action on the Fukaya category of a
surface; this is presumably well-known in certain circles. In contrast
with the Fukaya category, the triangulated categories constructed in
this paper are purely algebraic, and have finiteness properties
which are not apparent for the Fukaya category. There is, however, a
direct relation between the constructions in this paper and a variant
of the Fukaya category; see~\cite{AurouxBordered}.

This paper is structured as follows. In Section~\ref{sec:algebras} we define the algebras $\Blg(F)$;
these are more general than
the algebras from~\cite{LOT1}, since we allow $F$ to have more than one boundary
component, but are special cases of definitions
from~\cite{Zarev09:BorSut}. In Section~\ref{sec:bimodules} we define
the bimodules. In Section~\ref{sec:faithfulness} we prove faithfulness
of the action. In Section~\ref{sec:finite-generation} we discuss a
sense in which these categories are finitely generated, and the
decategorification of our action. We conclude, in
Section~\ref{sec:further}, with some further questions.

\subsection{Acknowledgements.} We thank T.~Cochran, E.~Grigsby and S.~Harvey for
helpful conversations, and C.~Clarkson and K.~Siegel for helpful comments on an
earlier version of this paper. The proof of faithfulness is inspired
by the argument in~\cite{KhS02:BraidGpAction}. We also thank the
Mathematical Sciences Research Institute and Columbia University for
hosting us during this
research. Finally, we thank the referee for many helpful comments.

\section{The algebras}\label{sec:algebras}

In the present section, we define two algebras associated to an arc
diagram~$\PMC$, denoted $\Blg(\PMC)$ (Section~\ref{sec:Blg}) and
$\Clg(\PMC)$ (Section~\ref{sec:Clg}). The second of these,
$\Clg(\PMC)$, is equipped with a differential, while the first,
$\Blg(\PMC)$, is not. These are both subalgebras of a more general
algebra $\Alg(\PMC)$, which we introduce in
Section~\ref{sec:Clg}. These algebras can be endowed with further
structure (notably, a kind of grading), which we will not need here;
see~\cite{LOT1}.

Before defining the algebras, we recall a convenient way of
representing surfaces.

\subsection{Arc diagrams}
Consider a connected, 
oriented surface $F$  of genus $g$ with $b>0$ boundary
components $Z_1,\dots,Z_b$, and suppose that each $Z_i$ is divided
into two closed arcs, $S_i^+$
and~$S_i^-$ (overlapping at their endpoints); write $S^+=\cup_i S_i^+$ and
$S^-=\cup_iS_i^-$. Choose a collection of pairwise-disjoint, embedded paths $\alpha_i$ in $F$
with $\bdy \alpha_i\subset S^+$ so that $F\setminus(\cup_i\alpha_i)$
is a union of disks, and the boundary of each disk contains
exactly one $S_i^-$. This implies that we have exactly $2(g+b-1)$ $\alpha$-curves. Place
a basepoint $z_i$ in each $S_i^-$.

Let $\{a_i,a'_i\}=\bdy\alpha_i$. We call 
\[
\PMC=(\overbrace{(Z_1,\dots,Z_b)}^Z,\overbrace{(\{a_1,a'_1\},\dots,\{a_n,a'_n\})}^{\Matching},\overbrace{(z_1,\dots,z_b)}^{\bbpt})
\]
an \emph{arc diagram} for $F$. Write $\CircPts=\{a_1,a'_1,\dots,a_{2(g+b-1)},a'_{2(g+b-1)}\}$. Here, the $Z_i$ are viewed as oriented circles.
For each $i$, the points $a_i$ and $a_i'$ are called a \emph{matched pair}. 

\begin{figure}
  \centering
  \includegraphics{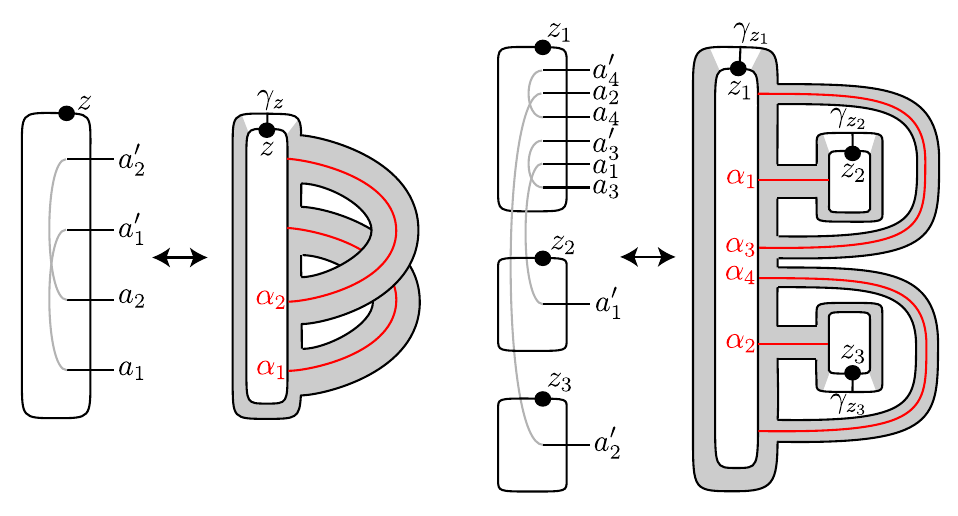}
  \caption{\textbf{Arc diagrams and their associated surfaces.} Left:
    an arc diagram specifying a once-punctured torus. Right: an arc
    diagram specifying the $3$-times punctured sphere. In each case, the subsurface
    $F(\PMC)\subset\PunctF(\PMC)$ is shaded.}
  \label{fig:arc-diagrams}
\end{figure}

From $\PMC$ we can build a standard model surface as follows. Thicken
the circles $Z_i$ in $Z$ to annuli $[0,1]\times Z_i$ and attach strips
($2$-dimensional $1$-handles) to each pair of points in $\Matching$
in the outer boundaries $\{1\}\times Z_i$ of the annuli. Call the result
$\PunctF(\PMC)$.  The basepoint $z_i$ in $Z_i$ gives an arc
$\gamma_{z_i}=[0,1]\times\{z_i\}\subset [0,1]\times Z_i$. Let
$F(\PMC)$ denote the result of cutting $\PunctF(\PMC)$ along the
$\gamma_{z_i}$. 
Let $S^+(\PMC)$ be the part of $\bdy F(\PMC)$ coming from $\bigcup_i\{0\}\times
Z_i$, together with the part corresponding to the $\gamma_{z_i}$,
and
let $S^-(\PMC)$ be the part of $\bdy F(\PMC)$ coming from
$\bigcup_i\{1\}\times Z_i$ (and the handles attached to it). 
See Figure~\ref{fig:arc-diagrams}.

The choice of the $\alpha_i$ identifies $F$ and $F(\PMC)$ canonically (up to isotopy).

\begin{remark}
  Let $\PMC=(Z,M,\bbpt)$ be an arc diagram for $F$. By definition,
  each circle in $Z$ contains one point
  $z_i\in\bbpt$. Moreover, performing surgery on $Z$ along the pairs of
  points in~$M$ gives a collection of circles each of which also contains a single
  $z_i$. Conversely, any triple $(Z,M,\bbpt)$ satisfying
  this condition comes from a surface.
\end{remark}

\begin{remark}
  We are considering a special case of Zarev's definition of arc
  diagrams \cite{Zarev09:BorSut}: he allows each $Z_i$ to be divided
  into $2n_i$ arcs for any $n_i\in\NN$.
\end{remark}

\subsection{The algebra \textalt{$\Blg(\PMC)$}{B(Z)}}\label{sec:Blg}
The algebras of interest are associated to arc diagrams. The algebra $\Blg(\PMC)$ has a basis over $\Field$ consisting of:
\begin{itemize}
\item One element $I_i$ for each pair of points $\{a_i,a'_i\}\in \Matching$.
\item One element $\rho$ for each nontrivial interval in each $Z_i\setminus \{z_i\}$ with endpoints in $\CircPts$. We will call these elements \emph{chords}. Given a chord $\rho$, let $\rho^-$ denote the initial point of $\rho$ (with respect to the orientation on $Z_i$), and let $\rho^+$ denote the terminal point of $\rho$.
\end{itemize}

The product on the algebra is given as follows:
\begin{itemize}
\item The $I_i$ are orthogonal idempotents, so $I_i^2=I_i$ and $I_iI_j=0$ if $i\neq j$.
\item $I_i\rho=\rho$ if $\rho^-$ is $a_i$ or $a'_i$; otherwise, $I_i\rho=0$. Similarly, $\rho I_j=\rho$ if $\rho^+$ is $a_j$ or $a'_j$; otherwise, $\rho I_j=0$.
\item For chords $\rho$ and $\sigma$, $\rho\sigma=0$ unless $\rho^+=\sigma^-$. If $\rho^+=\sigma^-$ then $\rho\sigma$ is the chord from $\rho^-$ to $\sigma^+$.
\end{itemize}

\begin{example}\label{eg:Alg-of-torus}
  There is a unique arc diagram $\PMC$ for the once-punctured torus,
  which is illustrated in Figure~\ref{fig:arc-diagrams}. The algebra
  $\Blg(\PMC)$ is $8$-dimensional, with basis
\[\{I_1,I_2,\rho_{1,2},\rho_{2,3},\rho_{3,4},\rho_{1,3},\rho_{2,4},\rho_{1,4}\}\]
  and multiplication table
\[
\begin{array}{l|@{\quad}llllllll}
  \hlx{v}
 \times & I_1 & I_2 & \rho_{1,2} & \rho_{2,3} & \rho_{3,4} & \rho_{1,3} & \rho_{2,4} & \rho_{1,4}\\
\hlx{hvv}
I_1 & I_1 & 0 & \rho_{1,2} & 0 & \rho_{3,4} & \rho_{1,3} & 0 & \rho_{1,4}\\
I_2 & 0 & I_2 & 0 & \rho_{2,3} & 0 & 0 & \rho_{2,4} & 0 \\
\rho_{1,2} & 0 & \rho_{1,2} & 0 & \rho_{1,3} & 0 & 0 & \rho_{1,4} & 0\\
\rho_{2,3} & \rho_{2,3} & 0 & 0 & 0 & \rho_{2,4} & 0 & 0 & 0\\
\rho_{3,4} & 0 & \rho_{3,4} & 0 & 0 & 0 & 0 & 0 & 0\\
\rho_{1,3} & \rho_{1,3} & 0 & 0 & 0 & \rho_{1,4} & 0 & 0 & 0\\
\rho_{2,4} & 0 & \rho_{2,4} & 0 & 0 & 0 & 0 & 0 & 0 \\
\rho_{1,4} & 0 & \rho_{1,4} & 0 & 0 & 0 & 0 & 0 & 0 \\
  \hlx{v}
\end{array}
\]
(When reading this table, the third entry in the top row, e.g.,
means that $I_1\rho_{1,2}=\rho_{1,2}$.)

We can encode this algebra more succinctly as
\[
\mathcenter{
\begin{tikzpicture}
  \node at (0,0) (I1) {$I_1$};
  \node at (2,0) (I2) {$I_2$};
  \draw[->, bend left=30] (I1) to node[above]{$\rho_{1,2},\ \rho_{3,4}$} (I2);
  \draw[->, bend left=30] (I2) to node[below]{$\rho_{2,3}$} (I1);
\end{tikzpicture}}
\mathcenter{\Big/(\rho_{2,3}\rho_{1,2}=\rho_{3,4}\rho_{2,3}=0).}
\]
\end{example}

\begin{example}\label{eg:antipodal}
  Let $F$ be a genus $g$ surface with one boundary component. One arc diagram for $F$ is obtained as follows. Label $4g+1$ points on a circle $Z$, in order, by
  \[
  z, a_1, \dots, a_{2g}, a'_1, \dots, a'_{2g}.
  \]
  The algebra $\Blg(\PMC)$ associated to this arc diagram $\PMC$ has
  idempotents $I_1,\dots,I_{2g}$. For convenience,
  define $I_{i+2g}=I_i$.  Then $\Blg(\PMC)$ is generated over $\Field$
  by $I_1,\dots,I_{2g}$ and elements $\rho_{i,j}$ for $1\leq i<j\leq
  4g$, with the relations:
  \begin{align*}
    I_i\rho_{i,j}I_j&=\rho_{i,j},\\
    I_i\rho_{j,k} = \rho_{j,k}I_i &= 0 \qquad \text{(in cases not
      covered above),}\\
    \rho_{i,j}\rho_{k,l}&=
    \begin{cases}
      \rho_{i,l} & \text{if }j=k,\\
      0 & \text{otherwise.}
    \end{cases}
  \end{align*}
Graphically, this is:
\[
\mathcenter{
\begin{tikzpicture}
  \node[anchor=mid] at (0,0) (I1) {$I_1$};
  \node[anchor=mid] at (3.25,0) (I2) {$I_2$};
  \node[anchor=mid] at (6.5,0) (dots) {$\cdots\!$};
  \node[anchor=mid] at (9.75,0) (Ig) {$I_{2g}$};
  \draw[->, bend left=20] (I1) to node[above]{$\rho_{1,2}$, $\rho_{2g+1,2g+2}$} (I2);
  \draw[->, bend left=20] (I2) to node[above]{$\rho_{2,3}$, $\rho_{2g+2,2g+3}$} (dots);
  \draw[->, bend left=20] (dots) to node[above]{$\rho_{2g-1,2g}$, $\rho_{4g-1,4g}$} (Ig);
  \draw[->, bend left=20] (Ig) to node[below]{$\rho_{2g,2g+1}$} (I1);
\end{tikzpicture}}
\bigg/
\left(
\begin{array}{l}
\rho_{i,i+1}\rho_{2g+i,2g+i+1}=0\\
\rho_{2g+i,2g+i+1}\rho_{i,i+1}=0
\end{array}\right).
\]
See also~\cite{AGW:KSandBordered}, where this algebra is related to the algebras in~\cite{KhS02:BraidGpAction}.
\end{example}

\begin{example}\label{eg:Alg-of-punctured-sphere}
  There is an arc diagram for the complement of $k>0$ disks in $S^2$
  that generalizes Figure~\ref{fig:arc-diagrams} (right).
  On one circle $Z_1$, label $3k-2$ points by
  \[
  z_1,a_1,b_1,a'_1,a_2,b_2,a'_2,\dots, a_{k-1},b_{k-1},a'_{k-1}.
  \]
  On each remaining $Z_i$ ($i=2,\dots,k$) place two points $z_i$ and
  $b'_i$. This represents a relabeling from
  Figure~\ref{fig:arc-diagrams}; see Figure~\ref{fig:arc-diag-relab}.

  \begin{figure}
    \centering
    \includegraphics{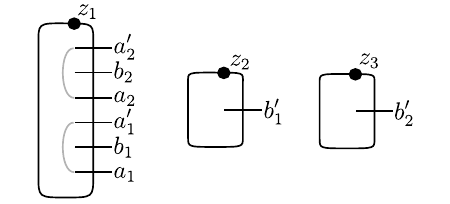}
    \caption{\textbf{Arc diagram for a thrice-punctured sphere.} This is a relabeling of the diagram from Figure~\ref{fig:arc-diagrams} (right).}
    \label{fig:arc-diag-relab}
  \end{figure}

  The associated algebra has idempotents $I_i$ ($i=1,\dots,k-1$) corresponding to the $\{a_i,a'_i\}$ and $J_i$ ($i=1,\dots,k-1$) corresponding to the $\{b_i,b'_i\}$. The algebra is given by
  \[
  \mathcenter{
  \begin{tikzpicture}
    \node[anchor=mid] at (0,0) (I1) {$I_1$};
    \node[anchor=mid] at (0,2) (J1) {$J_1$};
    \node[anchor=mid] at (3,0) (I2) {$I_2$};
    \node[anchor=mid] at (3,2) (J2) {$J_2$};
    \node[anchor=mid] at (5.25,0) (dots) {$\cdots\!$};
    \node[anchor=mid] at (7.5,0) (Inm1) {$I_{k-1}$};
    \node[anchor=mid] at (7.5,2) (Jnm1) {$J_{k-1}$};
    \draw[->, bend left=30] (I1) to node[left]{$\rho_{a_1,b_1}$} (J1);
    \draw[->, bend left=30] (J1) to node[right]{$\rho_{b_1,a'_1}$} (I1);
    \draw[->, bend right=20] (I1) to node[below]{$\rho_{a'_1,a_2}$} (I2);
    \draw[->, bend left=30] (I2) to node[left]{$\rho_{a_2,b_2}$} (J2);
    \draw[->, bend left=30] (J2) to node[right]{$\rho_{b_2,a'_2}$} (I2);
    \draw[->, bend right=20] (I2) to node[below]{$\rho_{a'_2,a_3}$} (dots);
    \draw[->, bend right=20] (dots) to node[below]{$\rho_{a'_{k-2},a_{k-1}}$} (Inm1);
    \draw[->, bend left=30] (Inm1) to node[left]{$\rho_{a_{k-1},b_{k-1}}$} (Jnm1);
    \draw[->, bend left=30] (Jnm1) to node[right]{$\rho_{b_{k-1},a'_{k-1}}$} (Inm1);
  \end{tikzpicture}}
\!\!\!\bigg/
\biggl(\begin{array}{l}
\rho_{b_i,a'_i}\rho_{a_i,b_i}=0\\
\rho_{a'_i,a_{i+1}}\rho_{a'_{i+1},a_{i+2}}=0
\end{array}\biggr).
\]
\end{example}

The following observation will be useful later:
\begin{lemma}\label{lem:opposite}
  Let $\PMC$ be an arc diagram and $-\PMC$ the arc diagram obtained by reversing the orientation of each circle $Z_i$ in $\PMC$. Then $\Blg(-\PMC)$ is the opposite algebra to $\Blg(\PMC)$.
\end{lemma}
\begin{proof}
  This is immediate from the definitions.
\end{proof}

\subsection{The algebra \textalt{$\Clg(\PMC)$}{C(Z)}}\label{sec:Clg}
Next we turn to the algebra $\Clg(\PMC)$. As mentioned in the
introduction, we give two different constructions, with $\Blg(\PMC)$
and with $\Clg(\PMC)$; either one gives a faithful action.  As such,
this section may be skipped at first reading.

Let $\PMC$ be an arc diagram
for a surface of genus $g$ with $b$ boundary components. Let
$n=2(g+b-1)$, so in particular the set $\CircPts$ of marked points has
$2n$ elements.

Consider $[0,1]\times (Z\setminus\bbpt)$. For each $i$ we can identify $[0,1]\times (Z_i\setminus z_i)$ with $[0,1]\times (-1,1)$. The points $\CircPts\cap Z_i$ give points $\{0\}\times \CircPts_i\subset \{0\}\times (-1,1)$ and $\{1\}\times \CircPts_i\subset \{1\}\times (-1,1)$.

A \emph{strand diagram} for $\PMC$ is a map $s\co \coprod_{i=1}^k
[0,1]\to [0,1]\times (Z\setminus\bbpt)$ (for some $k$), the components
of which we call strands, considered up
to reordering the strands, so that:
\begin{itemize}
\item $s$ maps $\coprod_{i=1}^k\{0\}$ to $\{0\}\times \CircPts\subset \{0\}\times Z$ and  $\coprod_{i=1}^k\{1\}$ to $\{1\}\times \CircPts\subset \{1\}\times Z$.
\item On each component of the source, $s$ is linear and has non-negative slope.
\item The map $s|_{\coprod_{i=1}^k\{0\}}$ is injective, as is the map $s|_{\coprod_{i=1}^k\{1\}}$.
\item For each matched pair $\{a_i,a'_i\}$, if there is a slope-zero
  strand (component of~$s$) starting at $(0,a_i)$ (respectively
  $(0,a'_i)$) then there is a slope-zero strand starting at $(0,a'_i)$
  (respectively $(0,a_i)$).
\item For each matched pair $\{a_i,a'_i\}$, if there is a positive-slope strand starting at $(0,a_i)$ (respectively $(0,a'_i)$) then there is no strand starting at $(0,a'_i)$ (respectively $(0,a_i)$).
\item  For each matched pair $\{a_i,a'_i\}$, if there is a positive-slope strand ending at $(1,a_i)$ (respectively $(1,a'_i)$) then there is no strand ending at $(1,a'_i)$ (respectively $(1,a_i)$).
\end{itemize}

\begin{figure}
  \centering
  \includegraphics[scale=.75]{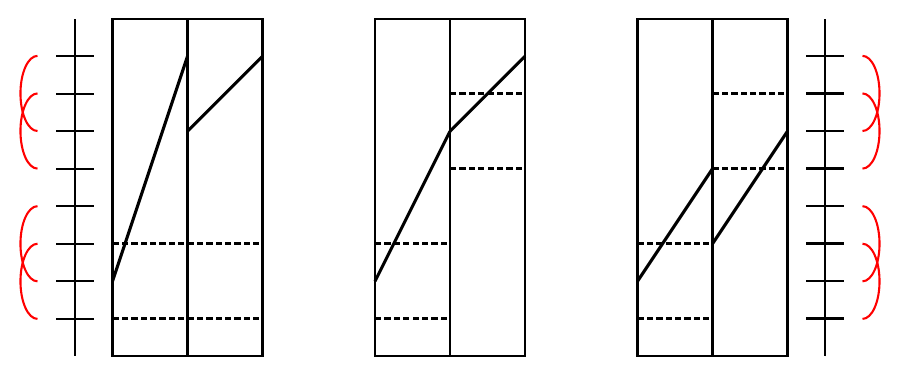}
  \caption{\textbf{Examples of $0$ products in $\Alg$.} The picture on
    the left illustrates the first two reasons the product can be
    zero, the picture in the middle illustrates the third and fourth
    reasons, and the picture on the right illustrates the last
    reason. When drawing elements of $\Alg$, we typically draw
    horizontal strands as dashed. This figure also appears in~\cite{LOT4}.}
  \label{fig:zero-prods}
\end{figure}

\begin{figure}
  \centering
  \includegraphics{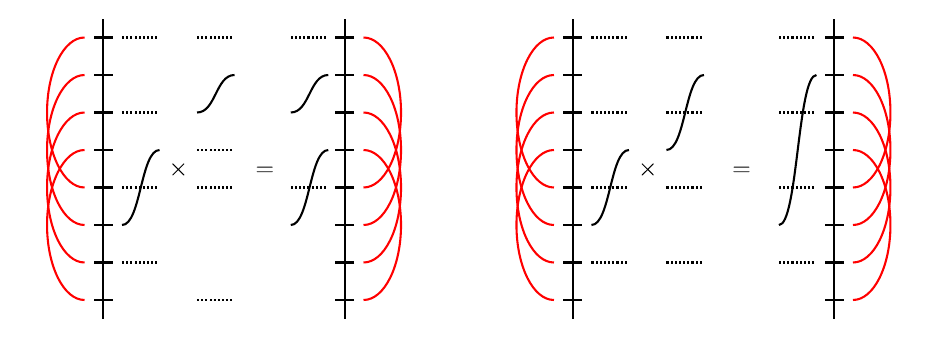}
  \caption{\textbf{Two nontrivial products.} 
    Both take place in $\Clg(\PMC)\subset
    \Alg(\PMC)$ for $\PMC$ the arc diagram from
    Example~\ref{eg:antipodal}. We have drawn the strands slightly
    curved, rather than straight, for artistic effect.}
  \label{fig:prod}
\end{figure}

Consider the $\Field$-vector space $\Alg(\PMC)$ generated by the strand diagrams. Define a product on this vector space as follows. Given $s, t\in\Alg(\PMC)$, the product of $s$ and $t$ is zero if
\begin{itemize}
\item there is a positive-slope strand in $s$ whose terminal endpoint is not the initial endpoint of a strand in $t$;
\item there is a positive-slope strand in $t$ whose initial endpoint is not the terminal endpoint of a strand in $s$;
\item there is a pair of slope-zero strands in $s$ neither of whose terminal endpoints is the initial endpoint of a strand in $t$;
\item there is a pair of slope-zero strands in $t$ neither of whose initial endpoints is the terminal endpoint of a strand in $s$;
 or
\item concatenating $s$ and $t$ end-to-end, there is a pair of
  piecewise-linear paths intersecting in two points (or
  equivalently, intersecting non-minimally).
\end{itemize}
See Figure~\ref{fig:zero-prods}.
In other cases, $s\cdot t$ is gotten by concatenating $s$ and $t$,
deleting any horizontal strands from $s$ (respectively $t$) which do
not match with strands in $t$ (respectively $s$), and pulling the
resulting piecewise-linear paths straight (fixing their
endpoints). See Figure~\ref{fig:prod}.


\begin{figure}
  \centering
  \includegraphics{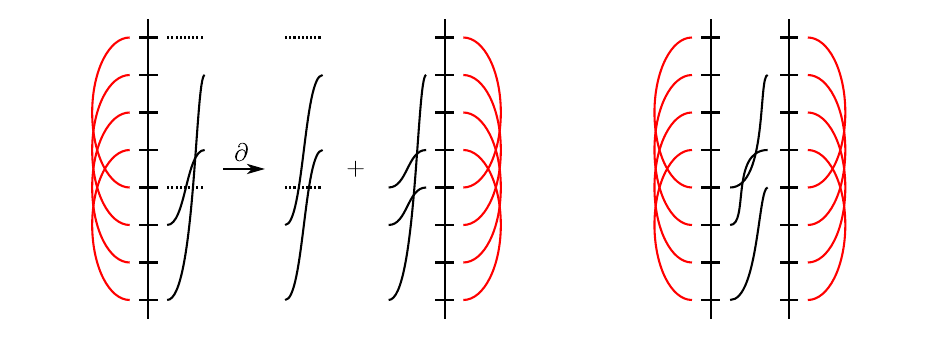}
  \caption{\textbf{A differential.} Left: an algebra element in $\Clg(\PMC)$, for $\PMC$ the arc diagram from Example~\ref{eg:antipodal}, and its differential. Right: a term which does not appear in the differential, because of a double-crossing. Again, we have drawn the strands slightly curved.}
  \label{fig:diff}
\end{figure}

Define a differential on $\Alg(\PMC)$ as follows. Given a strand diagram $s$ and a pair of intersecting strands $a,b$ in $s$, there is a unique (up to isotopy) way to resolve the intersection between $a$ and $b$ so that each resulting strand connects $\{0\}\times Z$ to $\{1\}\times Z$. If this resolution creates double-crossings between any pair of strands, let $s'_{a,b}=0$; otherwise, let $s'_{a,b}$ be the result of pulling straight the strands in the resolution and, if $a$ (respectively $b$) had slope $0$, deleting the slope-zero strand at $a'$ (respectively $b'$). Now, define
\[
\bdy(s)=\sum_{a,b\text{ intersect}}s'_{a,b}.
\]
See Figure~\ref{fig:diff}.

It is easy to verify that this multiplication and differential make $\Alg(\PMC)$ into a differential algebra.
The minimal idempotents for $\Clg(\PMC)$ are strand diagrams in which
all of the strands have slope~$0$, and so correspond to subsets of the
matched pairs in~$M$.

\begin{remark}
  It is easy to turn this geometric definition of $\Alg(\PMC)$ into a  combinatorial one; see, for instance,~\cite{LOT1}.
\end{remark}

The \emph{weight} of a strand diagram $s$ is the number of positive-slope strands in $s$ plus half the number of slope-zero strands in $s$. Let $\Alg(\PMC,k)$ be the subalgebra of $\Alg(\PMC)$ generated by strand diagrams of weight $k+n/2$. Then
\[
\Alg(\PMC)=\bigoplus_{k=-n/2}^{n/2}\Alg(\PMC,k).
\]

\begin{lemma}
  The algebra $\Blg(\PMC)$ is $\Alg(\PMC,-n/2+1).$
\end{lemma}
\begin{proof}
  This is immediate from the definitions.
\end{proof}

\begin{remark}
  The algebra $\Alg(\PMC,-n/2)$ is $\Field$ (generated by the empty
  strand diagram). The algebra $\Alg(\PMC,n/2)$ is quasi-isomorphic to
  $\Field$; compare Remark~\ref{rmk:KoszulDual}.
\end{remark}

\begin{definition}
  Let $\Clg(\PMC)=\Alg(\PMC,n/2-1).$
\end{definition}
In particular the algebra $\Clg(\PMC)$ has $n$ minimal idempotents, corresponding to the choices of $n-1$ of the $n$ matched pairs in $M$.

Given a chord $\rho$ in $\PMC$, let $c(\rho)\in\Clg(\PMC)$ be the sum of all ways of adding horizontal strands to $\rho$ to get an element of $\Clg(\PMC)$. (There are either $n-1$ such ways if the endpoints of $\rho$ are matched, or a single such choice if the endpoints of $\rho$ are not matched.)

\begin{example}
  For $\PMC$ the unique pointed matched circle for the torus, $\Clg(\PMC)\cong\Blg(\PMC)$, which is described explicitly in Example~\ref{eg:Alg-of-torus}.
\end{example}

\begin{example}\label{eg:planar-clg}
  Let $\PMC$ be the arc diagram from
  Example~\ref{eg:Alg-of-punctured-sphere} for the complement of $k$
  disks in $S^2$. The algebra $\Clg(\PMC)$ is quite large. However, as
  we will see, $\Clg(\PMC)$ is \emph{formal}; in fact, there is a map of
  algebras $f\co \Clg(\PMC)\to H_*(\Clg(\PMC))$ such that $f$ takes
  cycles to their homology classes. This means
  that in practice we can work with $H_*(\Clg(\PMC))$, which we
  describe explicitly below, instead of $\Clg(\PMC)$.

  To compute $H_*(\Clg(\PMC))$ (and see that $\Clg(\PMC)$ is formal),
  we use a little more terminology. Given a strand diagram
  $s\in\Clg(\PMC)$, the \emph{support} $\supp(s)$ of $s$ is the element of
  $H_1(Z,\CircPts)$ gotten by projecting $s$ to $Z$ and viewing the
  result as a $1$-chain. 

  As a first step towards understanding $H_*(\Clg(\PMC))$, let
  $M\subset \Clg(\PMC)$ be the $\Field$-subspace generated by strand
  diagrams $s$ such that $\supp(s)$ has multiplicity $>1$
  somewhere. Then $M$ is a differential ideal in $H_*(\Clg(\PMC))$.
  Further, $M$ is contractible, as in any arc
  diagram---see~\cite[Theorem 9]{LOT2}. So, it suffices to show
  that $\Clg'(\PMC)=\Clg(\PMC)/M$ is formal.

  Let $D$ be the $\Field$-subspace of $\Clg'(\PMC)$ generated by all
  strand diagrams $s\in\Clg'(\PMC)$ such that the interior of the
  $\supp(s)$ contains some
  point $b_i\in\CircPts$ which is occupied in the initial (and hence also
  in the terminal) idempotent. Then $D$ is a differential ideal in
  $\Clg(\PMC)$.  We claim that $D$ is contractible. To see this,
  consider a strand diagram $s\in D$. Such strand diagrams
  have one of two forms: either $s$ has some strand
  starting at a point $b_i$ (and hence also a strand ending at $b_i$)
  or it does not. The generators of the two types cancel in pairs:
  each generator of the first type occurs in the differential of a
  unique generator of the second type.

  Thus, we have reduced to considering
  $\Clg''(\PMC)=\Clg'(\PMC)/D$. It is now easy to see that the
  homology of $\Clg''(\PMC)$ is given by:
  \[
  \mathcenter{
  \begin{tikzpicture}
    \node[anchor=mid] at (0,0) (I1) {$I_1$};
    \node[anchor=mid] at (0,2) (J1) {$J_1$};
    \node[anchor=mid] at (3,0) (I2) {$I_2$};
    \node[anchor=mid] at (3,2) (J2) {$J_2$};
    \node[anchor=mid] at (5.25,0) (dots) {$\cdots\!$};
    \node[anchor=mid] at (7.5,0) (Inm1) {$I_{k-1}$};
    \node[anchor=mid] at (7.5,2) (Jnm1) {$J_{k-1}$};
    \draw[->, bend right=30] (J1) to node[left]{$\rho_{a_1,b_1}$} (I1);
    \draw[->, bend right=30] (I1) to node[right]{$\rho_{b_1,a'_1}$} (J1);
    \draw[->, bend left=20] (I2) to node[below]{$\rho_{a'_1,a_2}$} (I1);
    \draw[->, bend right=30] (J2) to node[left]{$\rho_{a_2,b_2}$} (I2);
    \draw[->, bend right=30] (I2) to node[right]{$\rho_{b_2,a'_2}$} (J2);
    \draw[->, bend left=20] (dots) to node[below]{$\rho_{a'_2,a_3}$} (I2);
    \draw[->, bend left=20] (Inm1) to node[below]{$\rho_{a'_{k-2},a_{k-1}}$} (dots);
    \draw[->, bend right=30] (Jnm1) to node[left]{$\rho_{a_{k-1},b_{k-1}}$} (Inm1);
    \draw[->, bend right=30] (Inm1) to node[right]{$\rho_{b_{k-1},a'_{k-1}}$} (Jnm1);
  \end{tikzpicture}}
\!\!\!\Bigg/
\Biggl(\begin{array}{l}
  \rho_{a_i',a_{i+1}}\rho_{b_i,a'_i}=0\\
  \rho_{a_{i},b_{i}}\rho_{a'_{i-1},a_{i}}=0\\
  \rho_{b_i,a'_i}\rho_{a_i,b_i}=0
\end{array}\Biggr).
\]
  (Here, $I_i$ corresponds to $\{a_i,a'_i\}$ not occupied and $J_i$
  corresponds to $\{b_i,b'_i\}$ not occupied. Each $\rho_{i,j}$ in
  the diagram actually stands for the homology classes of $\rho_{i,j}$
  in $H_*(\Clg(\PMC))$.)
  Further, the map $\Clg(\PMC)\to H_*(\Clg(\PMC))$ sending strand
  diagrams appearing in this homology to themselves and all other
  strand diagrams to $0$ is a map of algebras.
\end{example}

The following generalization of Lemma~\ref{lem:opposite} will be used
implicitly below:
\begin{lemma}\label{lem:opposite-gen}
  Let $\PMC$ be an arc diagram and $-\PMC$ the arc diagram obtained by
  reversing the orientation of each circle $Z_i$ in $\PMC$. Then
  $\Alg(-\PMC,i)$ is the opposite algebra to $\Alg(\PMC,i)$.
\end{lemma}
\begin{proof}
  This is immediate from the definitions.
\end{proof}

\begin{remark}
  \label{rmk:KoszulDual}
  It is not a coincidence that the algebra $\Clg(\PMC)$ from
  Example~\ref{eg:planar-clg} is formal: it follows from~\cite[Theorem
  9]{LOTHomPair} that $\Clg(\PMC)$ is always quasi-isomorphic to
  $\Blg(\PMC')$, where $\PMC'$ denotes the dual arc diagram to $\PMC$,
  as defined in Section~\ref{sec:diagrams}. (The reader may also
  notice a similarity between $\Blg(\PMC)$ and $\Blg(\PMC')$; it
  follows from results in~\cite{LOTHomPair} that these algebras are
  \emph{Koszul dual}; see also Remark~\ref{rmk:DDHalfId-is-dualizer}.)
\end{remark}

\begin{remark}\label{remark:nilpotent}
  Computations in $\Blg(\PMC)$ and $\Clg(\PMC)$ tend to be finite. In
  particular, both $\Blg(\PMC)$ and $\Clg(\PMC)$ are
  finite-dimensional. If we grade $\Blg(\PMC)$ and $\Clg(\PMC)$ by the
  total length (support) of an element, then all non-idempotent basic
  generators have positive grading. Thus, there is a number $N$,
  depending on $\PMC$, so that for any non-idempotent basic generators
  $a_1,\dots, a_N$ in $\Blg(\PMC)$ (respectively in $\Clg(\PMC)$),
  we have $a_1\cdots a_N=0$.
\end{remark}

\section{The bimodules}\label{sec:bimodules}
Let $\MCG_0(F)$ denote the mapping class group of $F$ fixing the boundary of $F$ pointwise. Our goal is to associate a bimodule $\CFDAa(\phi)$ to each element $\phi\in\MCG_0(F)$.
The definitions of the bimodules $\CFDAa(\phi)$ in~\cite{LOT2} and~\cite{Zarev09:BorSut}, even in the special case of interest to this paper, use holomorphic curves in a high symmetric product of a Riemann surface. We can work instead in the first symmetric product, making the whole story combinatorial, by taking advantage of a duality discussed in~\cite{LOTHomPair}. (In fact, there are two ways to do so, corresponding to using type \DD\ or type \AAm\ modules; we explain these in Sections~\ref{sec:D} and~\ref{sec:A}, respectively.)

\subsection{Diagrams for elements of the mapping class group}\label{sec:diagrams}

\begin{figure}
  \centering
  \includegraphics{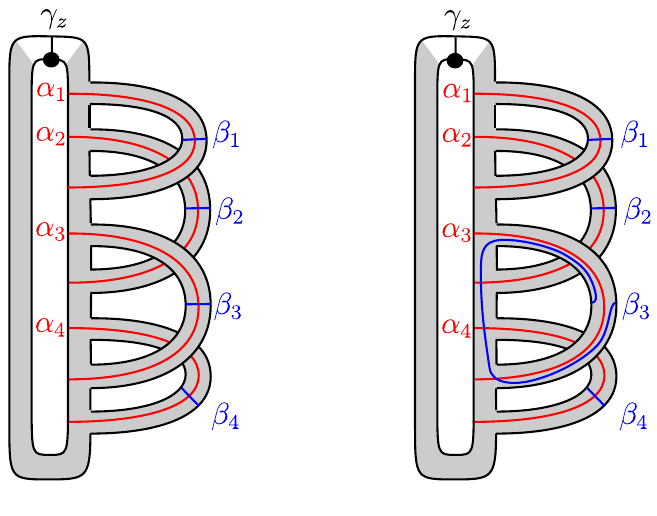}
  \caption{\textbf{Diagrams for mapping classes.} Left: a diagram for
    the identity map of the linear pointed matched circle. Right: a
    diagram for a (particular) Dehn twist. In each case, the subsurface
    $F(\PMC)\subset\PunctF(\PMC)$ is shaded.}
  \label{fig:mcg-diagrams}
\end{figure}

Fix an arc diagram $\PMC$, with $n$ pairs of matched points. As
discussed in Section~\ref{sec:algebras}, $\PMC$~specifies a surface
with boundary $\PunctF(\PMC)$ and a collection of arcs $\alpha_i$ in
$\PunctF(\PMC)$, whose complement is a union of disks. There is a dual
set of curves $\eta_i$ in $\PunctF(\PMC)$ so that 
\begin{itemize}
\item $\eta_i$ is contained in the handle of $\PunctF(\PMC)$ corresponding to $\alpha_i$ and
\item $\eta_i$ intersects $\alpha_i$ in a single point.
\end{itemize}
(See Figure~\ref{fig:mcg-diagrams}.)
Notice that $\{\eta_i\cap \bdy F(\PMC)\}$ is another arc diagram; we will call this the \emph{dual arc diagram to $\PMC$} and denote it $\PMC'$.

\begin{lemma}\label{lem:char-dual-curves}
  Up to isotopy, the $\eta_i$ are the unique curves in $F(\PMC)$ with
  boundary on $S^-$ and such that $\eta_i$ intersects $\alpha_i$ once
  and is disjoint from $\alpha_j$ for $i\neq j$.
\end{lemma}
\begin{proof}
  By definition, cutting along the $\alpha_i$ gives a disjoint union
  of disks. The boundary of each resulting disk will be divided into
  arcs coming from the original $S^-$ boundary, the original $S^+$
  boundary, and from the cut-open $\alpha$-curves.  The conditions on
  a pointed matched circle guarantee that there is a unique $S^-$
  interval on the boundary of each disk, so the $S^+$ and $\alpha$
  intervals necessarily alternate with each other.  The image of the
  $\eta_i$ in these disks are arcs in the interior that meet $S^-$ and
  one of the
  $\alpha$-curves.  These
  are uniquely characterized, up to isotopy. The result follows.
\end{proof}

Given $\phi\in\MCG_0(F(\PMC))$, viewing $F(\PMC)$ as a subsurface of $\PunctF(\PMC)$, we can act by $\phi$ on the $\eta$-curves, giving a new set of curves $\beta_i$. Write $\alphas=\alpha_1\cup\cdots\cup\alpha_n$ and $\betas=\beta_1\cup\cdots\cup\beta_n$. Let $\Diagram(\phi)=(\PunctF(\PMC),\alphas,\betas)$. Again, see Figure~\ref{fig:mcg-diagrams}. We will always assume that $\alphas\pitchfork\betas$; this is easy to arrange by deforming $\phi$ or $\{\alpha_i\}$ slightly.

The definitions of the bimodules will involve polygons in
$\Diagram(\phi)$. Assume, for convenience, that all of the
intersections between $\alphas$ and $\betas$ are right angles. Let 
\begin{align*}
  \bD^2&=\{x +iy\in\CC\mid x\geq 0, x^2+y^2\leq 1\}\\
  \gamma_R&=\bdy\bD^2\cap\{x+iy\in\CC\mid x\geq0\}\\
  \gamma_L&=\bdy\bD^2\cap\{x+iy\in\CC\mid x=0\}.
\end{align*}
Orient $\gamma_R$ and $\gamma_L$ from $-i$ to~$i$.

\begin{definition}
 Given chords $\rho_1,\dots,\rho_n$ in $\PMC$ and $\sigma_1,\dots,\sigma_m$ in $-\PMC'$, and points $x,y\in\alphas\cap\betas$, a \emph{polygon in $\Diagram(\phi)$ connecting $x$ to $y$ through $(\rho_1,\dots,\rho_n)$ and $(\sigma_1,\dots,\sigma_m)$} is a map $u\co \bD^2\to \Diagram(\phi)$ such that:
\begin{itemize}
\item $u(\gamma_L)\subset (\betas\cup\bdy\Diagram(\phi))$ and $u(\gamma_R)\subset (\alphas\cup\bdy\Diagram(\phi))$.
\item There are points $p_1,\dots,p_{2n}\in \gamma_R$ (respectively
  $q_1,\dots,q_{2m}\in \gamma_L$), appearing in that order as one traverses
  $\gamma_R$ (respectively~$\gamma_L$) from $-i$ to $i$, so that $u$ is an
  orientation-preserving immersion on
  $\bD^2\setminus\{p_1,\dots,p_{2n},q_1,\dots,q_{2m}\}$.  In particular,
  the image must (locally) make a right angle at $u(i)$ and $u(-i)$.
\item $u(-i)=x$ and $u(i)=y$.
\item For each $i$, 
$u([p_{2i+1},p_{2i+2}])=\rho_i$ and  $u([q_{2i+1},q_{2i+2}])=\sigma_i$; and except for these intervals, $u$ maps to the interior of $\Diagram(\phi)$.
\end{itemize}
\end{definition}

\begin{figure}
  \centering
  \includegraphics{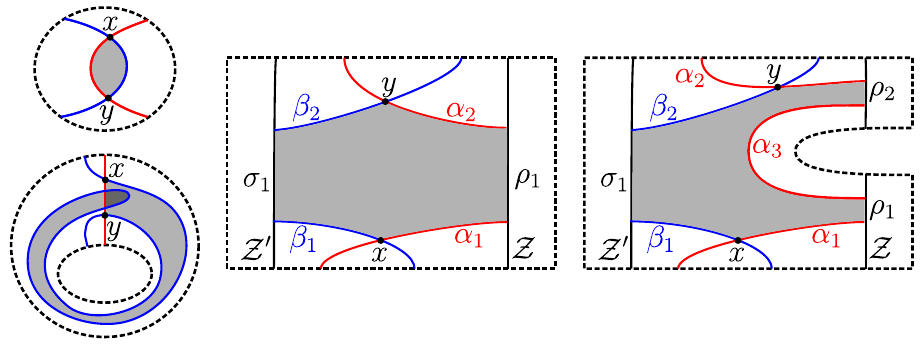}
  \caption{\textbf{Some polygons.} These are local pictures, that is,
  only a part of the diagram is shown (the boundary of which is
  indicated with dashed lines). Left: two polygons contributing to
  $n(x,y,(),())$, i.e., \emph{bigons}. In the lower of the two, the
  darker region is covered with multiplicity two. Center: a polygon contributing
to $n(x,y,(\rho_1),(\sigma_1))$. Right: a polygon contributing to
$n(x,y,\allowbreak (\rho_1,\rho_2),(\sigma_1))$.}
  \label{fig:polygons}
\end{figure}

See Figure~\ref{fig:polygons} for some sample polygons. Note that the sequence $(\sigma_1,\dots,\sigma_m)$ or $(\rho_1,\dots,\rho_n)$ (or both) may be empty. If both sequences are empty, we are counting the number of bigons between $x$ and $y$.

Call polygons $u$ and $v$ (connecting $x$ to $y$ and through
$(\rho_1,\dots,\rho_n)$ and $(\sigma_1,\dots,\sigma_m)$) \emph{equivalent} if
there is a diffeomorphism $w\co \bD^2\to\bD^2$ so that $v=u\circ
w$. 
Let $n(x,\allowbreak y,\allowbreak
(\sigma_1,\dots,\sigma_m),\allowbreak (\rho_1,\dots,\rho_n))\in\ZZ/2$
denote the number of equivalence classes of polygons connecting $x$ to
$y$ and through $(\rho_1,\dots,\rho_n)$ and $(\sigma_1,\dots,\sigma_m)$. (It is straightforward to check that the number of such polygons is always finite.)

\begin{lemma}\label{lem:polygon-is-polygon}
  Let $\cM(x,y,(\sigma_1,\dots,\sigma_m),(\rho_1,\dots,\rho_n))$ be
  the moduli space of pseudoholomorphic curves as in~\cite{LOT2}
  and~\cite{LOTHomPair} connecting $x$ to $y$ with asymptotics
  specified by the sequences
  $(\sigma_1,\dots,\sigma_m),(\rho_1,\dots,\rho_n)$. If
  this moduli space is
  $0$-dimensional then
\[
n(x,y,(\sigma_1,\dots,\sigma_m),(\rho_1,\dots,\rho_n))\equiv\#\cM(x,y,(\sigma_1,\dots,\sigma_m),(\rho_1,\dots,\rho_n))\pmod{2}.
\]
Otherwise, $n(x,y,(\sigma_1,\dots,\sigma_m),(\rho_1,\dots,\rho_n))=0$.
\end{lemma}
\begin{proof}
  This follows from the definitions and the Riemann mapping theorem.
\end{proof}

\begin{remark}
  Counting immersed polygons in $\Diagram(\phi)$ is combinatorial, and boils down to the combinatorics of gluing together components of $\Diagram(\phi)\setminus (\alphas\cup\betas)$.
\end{remark}

\subsection{Type \textalt{$D$}{D} modules}\label{sec:D}
The goal of this section is to associate a (differential) $\Clg(\PMC)$-bimodule $\ADM(\phi)$ to a strongly based mapping class $\phi\co F(\PMC)\to F(\PMC)$. We first define a $\Clg(\PMC')$-$\Clg(\PMC)$-bimodule $\DM(\phi)$ associated to $\phi$, and then define $\ADM(\phi)$ in terms of $\DM(\phi)$.

Given a point $x\in\alphas\cap\betas$ define $I_D(x)$ to be the
idempotent in $\Clg(\PMC')$ corresponding to the $\beta$-curves not
occupied by $x$, and $J_D(x)$ to be the idempotent in $\Clg(\PMC)$
corresponding to the $\alpha$-curves not occupied by $x$. Let
\[
\DM(\phi)=\bigoplus_{x\in\alphas\cap\betas}\Clg(\PMC')I_D(x)\otimes_\Field J_D(x)\Clg(\PMC).
\]
This is a $\Clg(\PMC')$-$\Clg(\PMC)$-bimodule. Abusing notation
imperceptibly, we let $x$ denote the generator for $\DM(\phi)$ corresponding to the intersection point $x$. Define a differential on $\DM(\phi)$ by
\[
\bdy(x)=\sum_{y\in\alphas\cap\betas}\sum_{\begin{subarray}{c}(\rho_1,\dots,\rho_n),\\(\sigma_1,\dots,\sigma_m)\end{subarray}}n(x,y,(\sigma_1,\dots,\sigma_m),(\rho_1,\dots,\rho_n))c(\sigma_1)\cdot\cdots\cdot c(\sigma_m)\cdot y\cdot c(\rho_n)\cdot\cdots\cdot c(\rho_1),
\]
and extending via the Leibniz rule
$\bdy(axb)=(\bdy(a))xb+a(\bdy(x))b+ax(\bdy(b))$.

\begin{lemma}\label{lem:DM-agrees}
  The bimodule $\DM(\phi)$ agrees with the bimodule
  $\CFDDa(\Diagram(\phi), n/2-1)$ as defined
  in~\cite{LOT2},~\cite{Zarev09:BorSut} and~\cite{LOTHomPair}.
\end{lemma}
\begin{proof}
  This is straightforward from the definitions and
  Lemma~\ref{lem:polygon-is-polygon}. (Note that in~\cite{LOTHomPair}
  we would have thought of $\DM(\phi)$ as a left module over
  $\Clg(-\PMC)$ and a right module over $\Clg(-\PMC')$; using the fact
  that $\Clg(-\PMC)\cong \Clg(\PMC)^\op$ we are viewing $\DM(\phi)$ as
  a right module over $\Clg(\PMC)$ and left module over $\Clg(\PMC')$.)
\end{proof}

\begin{proposition}
  If $\phi$ is isotopic relative to the boundary of $F(\PMC)$ to $\psi$ then $\DM(\phi)$ is homotopy equivalent to $\DM(\psi)$.
\end{proposition}
\begin{proof}
  This follows from the
  identification $\DM(\phi)\cong \CFDDa(\Diagram(\phi),n/2-1)$
  (Lemma~\ref{lem:DM-agrees}) and the corresponding invariance
  property of $\CFDDa(\Diagram(\phi),\allowbreak n/2-1)$. (It should
  also be possible to give a direct proof, since all of the objects
  involved are topological.)
\end{proof}

The bimodules $\DM(\phi)$ are not the ones promised in the introduction; indeed, they are bimodules over two different algebras. We perform one further algebraic operation to them. Let $\Id_{\PMC}$ denote the identity map of $F(\PMC)$. Then for $\phi\in\MCG_0(F(\PMC))$ define
\[
\ADM(\phi)=\Mor_{\Clg(\PMC')}(\DM(\Id_{\PMC}),\DM(\phi)),
\]
where $\Mor$ denotes the chain complex of left module maps
$\DM(\Id_{\PMC})\to\DM(\phi)$, which is a
$\Clg(\PMC)$-bimodule.\footnote{That is, $\Mor_{\Clg}(M,N)$ is
  generated by maps from $M$ to $N$ which respect the left module
  structure but not the right module structure or differential. The
  differential of such a map $f$ is given by
  $d(f)(x)=\bdy(f(x))+f(\bdy(x))$. The right action on
  $\Mor_{\Clg}(M,N)$ is given by $(f\cdot b)(x)=f(x)\cdot b$. The left
  action on $\Mor_{\Clg}(M,N)$ is given by $(b\cdot f)(x)= f(x\cdot
  b)$.}

\begin{proposition}\label{prop:bimodules-agree-D}
  The bimodule $\ADM(\phi)$ defined above agrees with the bimodule $\CFDAa(\phi,n/2-1)$ defined in~\cite{LOT2} (or~\cite{Zarev09:BorSut}).
\end{proposition}
\begin{proof}
The diagram $\Diagram(\phi)$ is an $\alpha$-$\beta$-bordered Heegaard diagram.
as in~\cite{LOTHomPair}. On the other hand,
$\Diagram(\Id_{-\PMC})\cup_{\PMC'}\Diagram(\phi)$ is an
$\alpha$-$\alpha$-bordered Heegaard diagram for $\phi$, in the sense
of~\cite{LOT2}. Thus, the pairing theorem for bordered Floer
homology expresses the bordered invariant $\CFDAa(\phi)$ as the
$\Ainf$ tensor product
  \begin{equation}
    \CFDAa(\phi)\simeq \CFAAa(\Diagram(\Id_{-\PMC}))\DTP_{\Alg(\PMC')}\CFDDa(\Diagram(\phi)).\label{eq:equal-1}
  \end{equation}
  (See~\cite[Section 7]{LOT2} for the pairing theorem and,
  for instance,~\cite{AinftyAlg} for a discussion of the
  $\Ainf$ tensor product.)
  In particular, taking $\phi=\Id$, the bimodules
  $\CFDDa(\Diagram(\Id_{\PMC}))$ and $\CFAAa(\Diagram(\Id_{-\PMC}))$
  are quasi-inverses to each other (in the sense of~\cite[Definition~\ref*{LOT2:def:QuasiInvertible}]{LOT2}. So, 
  \begin{multline}
      \Mor_{\Alg(\PMC')}(\CFDDa(\Diagram(\Id_{\PMC})),\Alg(\PMC'))\\
    \begin{aligned}[t]
       &\simeq \Mor_{\Alg(\PMC')}(\CFAAa(\Diagram(\Id_{-\PMC})) \DTP
       \CFDDa(\Diagram(\Id_{\PMC})),\CFAAa(\Diagram(\Id_{-\PMC})))\\
       &\simeq \CFAAa(\Diagram(\Id_{-\PMC})).
    \end{aligned}
    \label{eq:equal-2}
  \end{multline}
  Combining Equations~(\ref{eq:equal-1}) and~(\ref{eq:equal-2}) and
  using the fact that $\CFDDa(\Diagram(\phi),n/2-1)\cong \DM(\phi)$
  gives the result.
\end{proof}

\begin{corollary}\label{cor:gives-action-D}
  The bimodules $\ADM(\phi)$ satisfy the properties that
  $\ADM(\phi)\otimes\ADM(\psi)\simeq \ADM(\psi\circ \phi)$ and
  $\ADM(\Id)\simeq \lsub{\Clg(\PMC)}\Clg(\PMC)_{\Clg(\PMC)}$. In
  particular, they give an action of
  $\MCG_0(F(\PMC))$ on the derived category of right differential modules over $\Clg(\PMC)$.
\end{corollary}
\begin{proof}
  This follows from Proposition~\ref{prop:bimodules-agree-D} and the corresponding facts for $\CFDAa(\phi)$, which are proved in~\cite{LOT2}.
\end{proof}

\begin{remark}
  The bimodules $\ADM(\phi)$ are left- and right-projective
  (compare~\cite[Corollary~\ref*{LOT2:cor:tensor-projective}]{LOT2}),
  so the ordinary tensor
  product in
  Corollary~\ref{cor:gives-action-D} agrees with the derived tensor product.
\end{remark}

\begin{figure}
  \centering
  \includegraphics{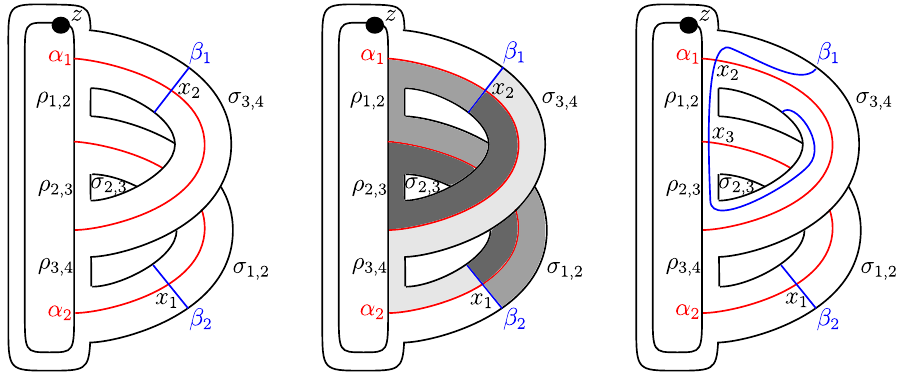}
  \caption{\textbf{The identity map of the torus and a Dehn twist.}
    Left: a diagram for the identity map of the torus. Center: three
    different polygons (rectangles) in the diagram, drawn with three
    different shadings. Right: a diagram for a Dehn twist around a particular
    essential curve in the torus.}
  \label{fig:torus-diags}
\end{figure}

\begin{example}\label{eg:D-id}
  Figure~\ref{fig:torus-diags} (left) shows a diagram for the identity map of
  the torus. There are three polygons (shown with
  different shadings in Figure~\ref{fig:torus-diags} in the middle), contributing 
  \begin{align*}
    \bdy x_1&=\sigma_{2,3}x_2\rho_{2,3}\\
    \bdy x_2&=\sigma_{1,2}x_1\rho_{1,2}+\sigma_{3,4}x_1\rho_{3,4}
  \end{align*}
  to the differential on $\DM(\Id)$.  No other polygons contribute to
  the differential: for polygons to contribute in this case, their
  boundaries must contain a single connected segment in each of $\PMC$
  and $\PMC'$, with multiplicity one. Any polygon whose image is the
  union of two components of
  $\Diagram(\Id)\setminus(\alphas\cup\betas)$ cannot contribute for
  idempotent reasons. The union of all three regions in
  $\Diagram(\Id)\setminus(\alphas\cup\betas)$ is represented by two
  different polygons, one contributing to
  $n(x_2, x_1, (\sigma_{1,4}), (\rho_{3,4},\rho_{2,3},\rho_{1,2}))$
  and one contributing to
  $n(x_2, x_1,\allowbreak
  (\sigma_{1,2},\sigma_{2,3},\sigma_{3,4}),\allowbreak
  (\rho_{1,4}))$.
  These cancel algebraically (each contributes
  $\sigma_{1,4}x_1\rho_{1,4}$ to $\bdy x_2$).

  Next, to compute $\ADM(\Id)$, we consider
  $\Mor_{\Clg(\PMC')}(\DM(\Id),\DM(\Id))$. As a left
  $\Clg(\PMC')$-module, $\DM(\Id)$ is generated by
  \[
  S=\{
  x_1, x_1\rho_{1,2}, x_1\rho_{1,3}, x_1\rho_{1,4}, x_1\rho_{3,4},
  x_2, x_2\rho_{2,3}, x_2\rho_{2,4}\}.
  \]
  Let $x$ be any element of $S_1=\{x_1, x_1\rho_{1,2},
  x_1\rho_{1,3}, x_1\rho_{1,4}, x_1\rho_{3,4}\}$. Then there is an
  element in
  $\Mor_{\Clg(\PMC')}(\DM(\Id),\DM(\Id))$ sending $x$ to any element
  of 
  \begin{multline*}
  \{
  x_1, x_1\rho_{1,2},
  x_1\rho_{1,3}, x_1\rho_{1,4}, x_1\rho_{3,4},
  \sigma_{2,4}x_1, \sigma_{2,4}x_1\rho_{1,2},
  \sigma_{2,4}x_1\rho_{1,3}, \sigma_{2,4}x_1\rho_{1,4},\\
  \sigma_{2,4}x_1\rho_{3,4},
  \sigma_{2,3}x_2, \sigma_{2,3}x_2\rho_{2,3}, \sigma_{2,3}x_2\rho_{2,4}
  \}
  \end{multline*}
  and sending all other elements of $S$ to $0$.

  Similarly, for $y$ any element of $S_2=\{x_2, x_2\rho_{2,3},
  x_2\rho_{2,4}\}$ there is an element of
  $\Mor_{\Clg(\PMC')}(\DM(\Id),\allowbreak \DM(\Id))$ sending $y$ to any of 
  \begin{multline*}
  \{
  x_2, x_2\rho_{2,3}, x_2\rho_{2,4},
  \sigma_{1,3}x_2, \sigma_{1,3}x_2\rho_{2,3},
  \sigma_{1,3}x_2\rho_{2,4},\\
  \sigma_{1,2}x_1, \sigma_{1,4}x_1, \sigma_{3,4}x_1,
  \sigma_{1,2}x_1\rho_{1,2}, \sigma_{1,4}x_1\rho_{1,2}, \sigma_{3,4}x_1\rho_{1,2},
  \sigma_{1,2}x_1\rho_{1,3},   \sigma_{1,4}x_1\rho_{1,3},   \sigma_{3,4}x_1\rho_{1,3},\\
  \sigma_{1,2}x_1\rho_{1,4},   \sigma_{1,4}x_1\rho_{1,4},   \sigma_{3,4}x_1\rho_{1,4}, 
  \sigma_{1,2}x_1\rho_{3,4},   \sigma_{1,4}x_1\rho_{3,4},   \sigma_{3,4}x_1\rho_{3,4}  
  \}
  \end{multline*}
  and sending all other elements of $S$ to $0$.

  The next step in computing $\ADM(\Id)$ is to compute the
  differential and module structure on
  $\Mor_{\Clg(\PMC')}(\DM(\Id),\DM(\Id))$. This is cumbersome,
  although explicit. Some
  examples of this form can be found in~\cite[Section 7]{LOTHomPair}
  and~\cite[Section 8]{LOT4}. In Section~\ref{sec:practical} we will
  give a more practical way of working with one of our algebra actions.
\end{example}

 \begin{example}
   Figure~\ref{fig:torus-diags} on the right also shows a diagram for
   a particular
   Dehn twist of the torus. The associated module $\DM(\phi)$ has
   three generators, $x_1$, $x_2$ and $x_3$, with differentials
   \begin{align*}
     \bdy(x_1)&=\sigma_{2,3}x_3\\
     \bdy(x_2)&=\sigma_{3,4}x_1\rho_{3,4}+x_3\rho_{1,2}\\
     \bdy(x_3)&=\sigma_{1,2}x_1\rho_{1,3}+\sigma_{1,3}x_2\rho_{2,3}.
   \end{align*}
   Unlike Example~\ref{eg:D-id}, where $\bdy^2=0$ was forced by
   products in the algebra being zero, one of the cases of $\bdy^2=0$
   here involves cancellation:
   \[
   \bdy^2(x_3)=\bdy(\sigma_{1,2}x_1\rho_{1,3}+\sigma_{1,3}x_2\rho_{2,3})
   =\sigma_{1,2}\sigma_{2,3}x_3\rho_{1,3}+\sigma_{1,3}x_3\rho_{1,2}\rho_{2,3}=0.
   \]
 \end{example}

\subsection{Type \textalt{$A$}{A} modules}\label{sec:A}

Let $\AM(\phi)$ be the $\Field$-vector space generated by $\Gen(\phi)=\alphas\cap\betas$. 
We will make $\AM(\phi)$ into a ($\Ainf$-) bimodule over $\Blg(\PMC')$ and $\Blg(\PMC)$. To start, define a left action of $\Blg(\PMC')$ and a right action of $\Blg(\PMC)$ on $\AM(\phi)$ as follows. Given $x\in \Gen(\phi)$ and idempotents $I\in \Blg(\PMC')$, $J\in\Blg(\PMC)$ corresponding to arcs $\alpha_i$ and $\beta_j$ respectively, we have
\[
I\cdot x\cdot J=
\begin{cases}
  x & \text{if }x\in\alpha_i\cap \beta_j\\
  0 & \text{otherwise.}
\end{cases}
\]

Next, given a chord $\rho$ in $\Blg(\PMC)$, define
\[
x\cdot\rho=\sum_{y\in\Gen(\phi)}n(x,y,(),(\rho))y.
\]
Similarly, given a chord $\sigma$ in $\Blg(\PMC')$ and another point $y\in\Gen(\HD)$, define
\[
\sigma\cdot x=\sum_{y\in\Gen(\phi)}n(x,y,(\sigma),())y.
\]
We will also denote $\sigma\cdot x$ by $m_{1,1,0}(\sigma, x)$ and 
$x\cdot \rho$ by $m_{0,1,1}(x,\rho)$; the reason
will become clear presently.

\begin{example}
  In the diagram for a Dehn twist of the torus in
  Figure~\ref{fig:torus-diags}, $x_2\rho_{1,2}=x_3$.
\end{example}

\begin{figure}
  \centering
  \includegraphics{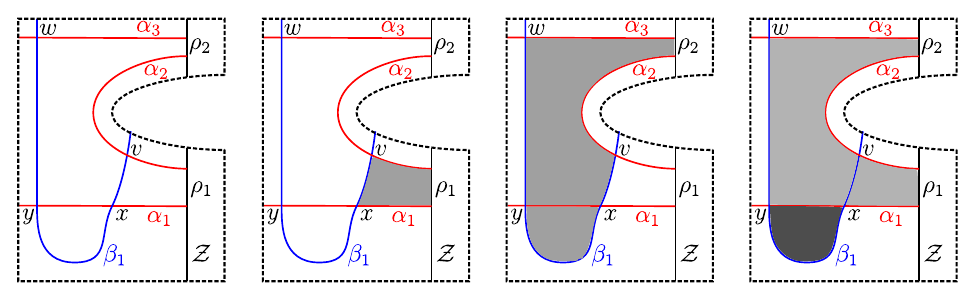}
  \caption{\textbf{Non-associativity of the action on $\AM(\phi)$.}
    This is a local example; only part of the diagram is drawn. There
    are products $x\cdot \rho_1=v$ and $v\cdot \rho_2=w$, given by the
  shaded regions in the second and third pictures,
  respectively. However, $x\cdot(\rho_1\rho_2)=0$: there is
  (obviously) no rectangle giving a nontrivial operation of this
  form. The resolution is that $\bdy(x)=y$ and
  $m_3(y,\rho_1,\rho_2)=w$; these operations are given by the darkly
  and lightly shaded regions in the fourth picture, respectively.}
  \label{fig:non-assoc}
\end{figure}

In general, the action we have defined so far may not be associative;
see Figure~\ref{fig:non-assoc}. As the notation suggests, we should
really think of $\AM(\phi)$ as an $\Ainf$-bimodule. As a warm-up,
define a differential on $\AM(\phi)$ by counting bigons:
\[
\bdy(x)=\sum_{y\in\Gen(\phi)}n(x,y,(),())y.
\]
It is straightforward to verify that $\bdy^2=0$, and the reader to
whom this is unfamiliar is encouraged to do so. We will also denote
$\bdy(x)$ as $m_{0,1,0}(x)$.

More generally, given a sequence of chords $\sigma_1,\dots,\sigma_n$ in $\Blg(\PMC')$ and $\rho_1,\dots,\rho_n$ in $\Blg(\PMC)$, and generators $x,y\in\Gen(\phi)$ define
\[
m_{m,1,n}(\sigma_m,\dots,\sigma_1,x,\rho_1,\dots,\rho_n)=\sum_{y\in\Gen(\phi)}n(x,y,(\sigma_1,\dots,\sigma_m),(\rho_1,\dots,\rho_n))y.
\]
Extend this multi-linearly to a map
\[
m_{m,1,n}\co
\overbrace{\Blg(\PMC')\otimes\cdots\otimes\Blg(\PMC')}^{m\text{ copies}}\otimes
\AM(\phi)\otimes \overbrace{\Blg(\PMC)\otimes\cdots\otimes\Blg(\PMC)}^{n\text{ copies}}\to \AM(\phi).
\]

\begin{lemma}
  These $m_{i,1,j}$ endow $\AM(\phi)$ with the structure of an $\Ainf$-bimodule.
\end{lemma}
\begin{proof}
  This is not too hard to prove combinatorially, but also follows from the analysis in~\cite{LOT1} and the Riemann mapping theorem. 
\end{proof}

\begin{remark}
  Even if $m_{0,1,0}=0$, so that $m_{1,1,0}$ and $m_{0,1,1}$ make $\AM(\phi)$ into an honest
  bimodule, there is a lot of additional information in the higher
  $\Ainf$-operations. See, for instance, Example~\ref{eg:A-id}.
\end{remark}

As with the bimodules $\DM(\phi)$ in Section~\ref{sec:D}, 
the bimodules $\AM(\phi)$ are not the ones promised in the introduction. Let $\Id_{\PMC}$ denote the identity map of $\PunctF(\PMC)$. Then for $\phi\in\MCG_0(F(\PMC))$ define
\begin{equation}
\DAM(\phi)=\Mor_{\Blg(\PMC')}(\AM(\Id_{\PMC}),\AM(\phi)),\label{eq:def-of-DAM}
\end{equation}
where $\Mor$ denotes the chain complex of left $\Ainf$-module
morphisms (whose cycles are the $\Ainf$-homomorphisms); see, for
instance,~\cite[Chapter 2]{LOT2}. Note that the right actions by
$\Blg(\PMC)$ on $\AM(\Id)$ and $\AM(\phi)$ give $\DAM(\phi)$ the
structure of an ($\Ainf$) $\Blg(\PMC)$-bimodule.

\begin{proposition}\label{prop:bimodules-agree-A}
  The bimodule $\DAM(\phi)$ defined above agrees with the bimodule
  $\CFDAa(\phi,\allowbreak -n/2+1)$ defined in~\cite{LOT2} (or~\cite{Zarev09:BorSut}).
\end{proposition}
\begin{proof}
  The proof is essentially the same as the proof of
  Proposition~\ref{prop:bimodules-agree-D}.
\end{proof}
\begin{corollary}
  The bimodules $\DAM(\phi)$ satisfy the properties that
  $\DAM(\phi)\DTP\DAM(\psi)\simeq \DAM(\psi\circ \phi)$ (where
  $\DTP$ denotes the $\Ainf$ tensor product) and $\DAM(\Id)\simeq \lsub{\Blg(\PMC)}\Blg(\PMC)_{\Blg(\PMC)}$. In particular, the bimodules $\DAM(\phi)$ give an action of $\MCG_0(F(\PMC))$ on the $\Ainf$-homotopy category of right $\Ainf$-modules over $\Blg(\PMC)$.
\end{corollary}
\begin{proof}
  Similarly to Corollary~\ref{cor:gives-action-D}, this follows from Proposition~\ref{prop:bimodules-agree-A} and the corresponding facts for $\CFDAa(\phi)$, which are proved in~\cite{LOT2}.
\end{proof}

\begin{example}\label{eg:A-id}
  For the identity map $\Id$ of the torus, using the diagram from
  Figure~\ref{fig:torus-diags} $\AM(\Id)$ has two generators $x_1$ and
  $x_2$. The differential and ordinary product are both trivial. There
  are, however, obvious higher products given by the rectangles in
  Figure~\ref{fig:torus-diags}, of the forms:
  \begin{align*}
    m_3(\sigma_{1,2},x_2,\rho_{1,2})&=x_1\\
    m_3(\sigma_{3,4},x_2,\rho_{3,4})&=x_1\\
    m_3(\sigma_{2,3},x_1,\rho_{2,3})&=x_2.
  \end{align*}
  This is not the end of the story; indeed, with only these products,
  $\AM(\Id)$ would not satisfy the $\Ainf$-relations. For instance,
  there is an operation
  \[
  m_4(\sigma_{2,3},\sigma_{1,2},x_2,\rho_{1,3})=x_2.
  \]
  To see this, consider the union of the regions abutting $\rho_{1,2}$
  and $\rho_{2,3}$. Make a cut in this region from $x_2$ along the
  $\beta$-curve to the boundary. The result is a polygon, from $x_2$
  to itself, through the chord $\rho_{1,3}$ on one boundary component
  and the chords $\sigma_{2,3}$ and
  $\sigma_{1,2}$ on the other boundary component. (This operation is
  also forced by the $\Ainf$-relation.)

  Similarly, there are higher products:
  \begin{align*}
    m_5(\sigma_{1,3},\sigma_{1,2},x_2,\rho_{1,3},\rho_{1,2})&=x_1\\
    m_6(\sigma_{2,3},\sigma_{1,3},\sigma_{1,2},x_2,\rho_{1,3},\rho_{1,3})&=x_2\\
    m_7(\sigma_{1,3},\sigma_{1,3},\sigma_{1,2},x_2,\rho_{1,3},\rho_{1,3},\rho_{1,2})&=x_1\\
    &\vdots\\
    m_5(\sigma_{3,4},\sigma_{2,3},\sigma_{1,2},x_2,\rho_{1,4})&=x_1\\
    m_6(\sigma_{2,4},\sigma_{2,3},\sigma_{1,2},x_2,\rho_{1,4},\rho_{2,3})&=x_2\\
    m_7(\sigma_{1,4},\sigma_{2,3},\sigma_{1,2},x_2,\rho_{1,4},\rho_{2,3},\rho_{1,2})&=x_1\\
    &\vdots,
  \end{align*}
  as well as several more infinite families, and 
  similar infinite families starting from $x_1$.
  
  It would be natural to compute $\DAM(\Id)$ next, via
  Equation~\eqref{eq:def-of-DAM}. This is tedious (and infinite); we
  will give a better method for computing $\DAM(\phi)$ from
  $\AM(\phi)$ in the next section. In particular, by
  Corollary~\ref{cor:enough-to-compute}, we have written down enough
  of $\AM(\phi)$ to characterize $\DAM(\phi)$ (as well as $\AM(\phi)$).
\end{example}

\subsection{Practical computations}\label{sec:practical}
As Example~\ref{eg:D-id} illustrates, computing the
bimodules $\ADM(\phi)$ and $\DAM(\phi)$ from the modules $\DM(\phi)$ and
$\AM(\phi)$ is quite cumbersome, and computing the tensor products
$\ADM(\phi_1)\otimes\ADM(\phi_2)$ or $\DAM(\phi_1)\DTP\DAM(\phi_2)$
would be even more so. In this section, we give a reformulation of the
bimodules $\DAM(\phi)$ which is better suited for computations. The
key tool is the type \DD\ bimodule associated to the diagram
$\Diagram(\Id)$ in the second to lowest $\SpinC$-structure. (The type
\DD\ bimodules we have worked with so far are in the second to  highest
$\SpinC$-structure.)

Call a chord in $\Blg(\PMC)$ \emph{short} if it connects adjacent
points in~$\CircPts$. Let $\ShortChords(\PMC)$ denote the set of
short chords in
$\PMC$. The diagram $\Diagram(\Id_\PMC)$ sets up a correspondence between
$\ShortChords(\PMC)$ and $\ShortChords(\PMC')$ as follows: two chords
correspond if they lie on the boundary of a single connected component
of $\Diagram(\Id_\PMC)\setminus (\alphas\cup\betas)$. Given a short chord
$\xi\in\ShortChords(\PMC)$ let $\xi'$ be the corresponding short chord
in $\ShortChords(\PMC')$.

\begin{definition}
  Given an arc diagram $\PMC$, let $\DDHalfId$ denote the
  $\Blg(\PMC)\Hyph\Blg(\PMC')$\hyp bimodule defined as
  follows. The bimodule $\DDHalfId$ has one generator $x_i$ for each
  matched pair $\{a_i,a'_i\}$ in $\Matching$. Let $I(x_i)$ be the
  idempotent in $\Blg(\PMC)$ corresponding to $\{a_i,a'_i\}$ and
  $J(x_i)$ the idempotent in $\Blg(\PMC')$ corresponding to
  $\{a_i,a'_i\}$. Let
  \[
  \DDHalfId=\bigoplus_i \Blg(\PMC)I(x_i)\otimes J(x_i)\Blg(\PMC').
  \]
  Abusing notation, we also let $x_i$ denote a generator of the summand
  corresponding to~$x_i$.
  Define a differential on $\DDHalfId$ by 
  \[
  \bdy(x_i)=\sum_{j}\sum_{\xi\in\ShortChords(\PMC)} I(x_i)\cdot \xi
  \cdot x_j\cdot \xi' \cdot J(x_i),
  \]
  and extending via the Leibniz rule. (Note that most terms in the sum
  defining $\bdy(x_i)$ vanish for idempotent reasons.)
\end{definition}

\begin{proposition}\label{prop:comp-DD-half-id}
  The bimodule $\DDHalfId$ is homotopy equivalent to
  $\CFDDa(\Diagram(\Id_{-\PMC}),\allowbreak -n/2+1)$.
\end{proposition}
\begin{proof}
  The identification of generators is given as follows: the generator
  $x_i$ for $\DDHalfId$ corresponding to the matched pair
  $\{a_i,a'_i\}$ corresponds to the generator $\x\subset
  \alphas\cap\betas$ for $\CFDDa(\Diagram(\Id_{-\PMC}))$ consisting of
  $\{\alpha_j\cap\beta_j\mid j\neq i\}$.  Each term in the
  differential on $\DDHalfId$ corresponds to an embedded hexagon in
  $\Diagram(\Id_{-\PMC})$, and hence corresponds to a term in the differential
  on $\CFDDa(\Diagram(\Id_{-\PMC}),-n/2+1)$. So, it remains to show that there
  are no other terms in the differential on
  $\CFDDa(\Diagram(\Id_{-\PMC}),-n/2+1)$. We will do this by showing that there
  are no other index $1$ positive domains whose boundaries in $\PMC$
  and~$\PMC'$ are such that they can contribute to the differential.

  Writing $\Diagram(\Id_{-\PMC})=(\PunctF(\PMC),\alphas,\betas)$, each
  component of $F(\PMC)\setminus(\alphas\cup\betas)$ is a hexagon,
  with two sides contained in $\alphas$, two sides contained in
  $\betas$, one side in $\PMC$ and one side in $\PMC'$. 

  For a domain $B$ to contribute $a'\cdot y\cdot a$ to $\bdy x$ we
  must have
  \[
  e(B)+n_\x(B)+n_\y(B)-\iota(a)-\iota(a')=-1.
  \]
  For generators $a\in\Blg(\PMC)$, we have $\iota(a)=0$ if $a$ is an
  idempotent and $-1/2$ if $a$ is not.  All non-trivial domains in
  $\Diagram(\Id_{-\PMC})$
  intersect both $\PMC$ and $\PMC'$, so
  $\iota(a)=\iota(a')=-1/2$. Thus, for a domain $B$ to contribute, it must
  have
  \[
  e(B)+n_\x(B)+n_\y(B)=0.
  \]

  Fix generators $\x=\{x_1,\dots,x_{2n-1}\}$ and
  $\y=\{y_1,\dots,y_{2n-1}\}$ for $\DDHalfId$. Reordering $\x$ and
  $\y$ if necessary, we may assume $x_i=y_i$ for $i<2n-1$. There are
  two cases: either $x_{2n-1}=y_{2n-1}$ or $x_{2n-1}\neq
  y_{2n-1}$. For simplicity, we will treat these two cases separately.
  
  \textbf{Case 1.} $x_{2n-1}=y_{2n-1}$. There is one point $p\in
  \alphas\cap\betas$ not appearing in $\x=\y$. Let $R_1,\dots,R_4$
  denote the four components of
  $\PunctF(\PMC)\setminus(\alphas\cup\betas)$ which have $p$ as a
  corner. Note that $R_i\neq R_j$ for $i \neq j$ unless $R_i$ contains
  a basepoint.

  If $R$ is some component of
  $\PunctF(\PMC)\setminus(\alphas\cup\betas)$ other than
  $R_1,\dots,R_4$ then
  \[
  e(R)+n_\x(R)+n_\y(R)=-1/2+1=1/2.
  \]
  By contrast, for the regions $R_1,\dots,R_4$,
  \[
  e(R_i)=-1/2+1/2=0.
  \]
  Thus, for any positive domain $B$, $e(B)+n_\x(B)+n_\y(B)\geq 0$, with
  equality if and only if $B$ is a linear combination of
  $R_1,\dots,R_4$. Thus, for
  grading reasons, the only domains that could contribute in this case
  are linear combinations of
  $R_1,\dots,R_4$.
  Because the algebra element on each side must be a
  single, connected chord, the multiplicity of each $R_i$ must be $0$
  or $1$. 
So, the rest of the argument boils down to the
  combinatorics of gluing together $\leq 4$ hexagons, each with two
  boundary components labeled
  $\alpha$, two labeled $\beta$, one labeled $\PMC$ and one labeled
  $\PMC'$, and gluing allowed along the $\alpha$- and $\beta$-boundary
  components.

  If we number the $R_i$ counter-clockwise around~$p$, say,
  with $R_1$ separated from $R_2$ by a $\beta$-arc
  then the only such linear combinations which could give domains in
  $\pi_2(\x,\y)$ are $R_1+R_2$, $R_2+R_3$, $R_3+R_4$, $R_4+R_1$ and
  $R_1+R_2+R_3+R_4$. If
  $R_1+R_2$ gives a domain then the chords in $\PMC$ corresponding to
  $R_1$ and $R_2$ must be consecutive. Such a domain has no
  holomorphic representative compatible with the idempotents, as in
  Example~\ref{eg:D-id}. The cases
  $R_2+R_3$, $R_3+R_4$ and $R_4+R_1$ are similar. For
  $R_1+R_2+R_3+R_4$, it is not possible for the domain to have
  connected boundary in $\PMC$ (or $\PMC'$).

  Thus, no domains from this case contribute to the differential on
  $\CFDDa(\Diagram(\Id_{-\PMC}))$.

  \textbf{Case 2.} $x_{2n-1}\neq y_{2n-1}$. As before,
  all regions have $e(R)+n_\x(R)+n_\y(R)\geq 0$. Moreover, equality
  only occurs for regions containing both $x_{2n-1}$ and $y_{2n-1}$ on
  their boundaries. There can be at most three such regions not
  containing basepoints; Figure~\ref{fig:torus-diags} (on the left) is
  the essentially unique case in which there are three. Let
  $R_1,R_2,R_3$ denote the three such regions (if three exist). The
  only linear combinations of $R_1,R_2,R_3$ giving domains in
  $\pi_2(\x,\y)$ with multiplicities $0$ or~$1$ everywhere in $\PMC$
  and $\PMC'$ are $R_1$, $R_2$, $R_3$, and $R_1+R_2+R_3$. The cases
  $R_1$, $R_2$ and $R_3$ contribute terms that occur in the
  differential on $\DDHalfId$. If $R_1+R_2+R_3$ exists then its
  geometry is exactly as in the genus~$1$ case
  (Figure~\ref{fig:torus-diags}). Thus, as in Example~\ref{eg:D-id},
  there are two cancelling holomorphic representatives.

  This concludes the proof of Proposition~\ref{prop:comp-DD-half-id}.
\end{proof}

\begin{corollary}
  The bimodule $\DDHalfId\DTP_{\Blg(\PMC')}\AM(\phi)$ is
  $\Ainf$-homotopy equivalent to $\DAM(\phi)$.
\end{corollary}
\begin{proof}
  This follows from Proposition~\ref{prop:comp-DD-half-id} and the
  definitions, similarly to the proof of
  Proposition~\ref{prop:bimodules-agree-A}.
\end{proof}

We have been using the notation $\DTP$ to denote the $\Ainf$ tensor
product. The resulting chain complexes are almost always
infinite-dimensional. For cases under consideration, however, there is
a smaller model for the $\Ainf$ tensor product, which we
denote~$\DT$. We refer the reader to~\cite{LOT2} for the definition.

\begin{example}\label{eg:DA-id}
  Continuing Example~\ref{eg:A-id}, we are now in a position to
  compute the bimodule $\DAM$ for the identity map of the torus. In this
  case, the bimodule $\DDHalfId$ has generators $w_1$ and $w_2$, with 
  \begin{align*}
    \bdy w_1&=\rho_{2,3}w_2\sigma_{2,3}\\
    \bdy w_2&=\rho_{1,2}w_1\sigma_{1,2}+\rho_{3,4}w_1\sigma_{3,4}.
  \end{align*}
  (This is, not coincidentally, the same as the bimodule $\DM(\Id)$
  from Example~\ref{eg:D-id}.)  Taking the $\DT$ tensor product with
  the bimodule $\AM(\Id)$ from Example~\ref{eg:A-id}, we get a bimodule
  with generators $w_1\otimes x_1$ and $w_2\otimes x_2$, and
  operations:
  \begin{align*}
    m_2(w_2\otimes x_2, \rho_{1,2})&=\rho_{1,2}w_1\otimes x_1 &
    m_2(w_2\otimes x_2, \rho_{3,4})&=\rho_{3,4}w_1\otimes x_1\\
    m_2(w_1\otimes x_1, \rho_{2,3})&=\rho_{2,3}w_2\otimes x_2 &
    m_2(w_2\otimes x_2, \rho_{1,3})&=\rho_{1,3}w_2\otimes x_2\\
    m_2(w_1\otimes x_1, \rho_{2,4})&=\rho_{2,4}w_1\otimes x_1 &
    m_2(w_2\otimes x_2, \rho_{1,4})&=\rho_{1,4}w_1\otimes x_1.
  \end{align*}
  This is exactly the $\Blg(\PMC)$-bimodule $\Blg(\PMC)$---as we
  expected.

  (The module is projectively generated on the left, like the type $D$
  cases above. On the right, it is an $\Ainf$-module, like the type
  $A$ cases.)

  We show a few examples of how these operations arise. The operation   
  $m_2(w_2\otimes x_2, \rho_{1,4})=\rho_{1,4}w_1\otimes x_1$ comes
  from a diagram of the following form:
  \[
  \begin{tikzpicture}
    \node[anchor=mid] at (0,0) (tl) {$w_2$};
    \node[anchor=mid] at (2,0) (tr) {$x_2$};
    \node at (4,0) (trr) {};
    \node at (0,-1) (delta) {$\bdy$};
    \node at (2,-2) (m) {$m_3$};
    \node[anchor=mid] at (0,-3) (bl) {$w_1$};
    \node[anchor=mid] at (2,-3) (br) {$x_1$};
    \node at (-2,-3) (bll) {};
    \draw[->] (tl) to (delta);
    \draw[->] (delta) to (bl);
    \draw[->] (tr) to (m);
    \draw[->] (m) to (br);
    \draw[->, bend right=15] (delta) to node[below, sloped]{\lab{\rho_{1,2}}} (bll);
    \draw[->, bend left=15] (delta) to node[below, sloped]{\lab{\sigma_{1,2}}} (m);
    \draw[->, bend right=10] (trr) to node[below, sloped]{\lab{\rho_{1,2}}} (m);
  \end{tikzpicture}
  \]
  using the higher product $m_3(\sigma_{1,2},x_2,\rho_{1,2})=x_1$.

  The operation   
  $m_2(w_2\otimes x_2, \rho_{1,2})=\rho_{1,2}w_1\otimes x_1$ comes
  from a diagram of the following form:
  \[
  \begin{tikzpicture}
    \node[anchor=mid] at (0,0) (tl) {$w_2$};
    \node[anchor=mid] at (4,0) (tr) {$x_2$};
    \node at (6,0) (trr) {};
    \node at (0,-1) (delta1) {$\bdy$};
    \node at (0,-2) (delta2) {$\bdy$};
    \node at (0,-3) (delta3) {$\bdy$};
    \node at (4,-4) (m) {$m_5$};
    \node at (-4,-4) (mult) {multiply};
    \node[anchor=mid] at (0,-5) (bl) {$w_1$};
    \node[anchor=mid] at (4,-5) (br) {$x_1$};
    \node at (-5.5,-5) (bll) {};
    \draw[->] (tl) to (delta1);
    \draw[->] (delta1) to (delta2);
    \draw[->] (delta2) to (delta3);
    \draw[->] (delta3) to (bl);
    \draw[->] (tr) to (m);
    \draw[->] (m) to (br);
    \draw[->, bend right=15] (delta1) to node[below, sloped]{\lab{\rho_{1,2}}} (mult);
    \draw[->, bend right=10] (delta2) to node[below, sloped]{\lab{\rho_{2,3}}} (mult);
    \draw[->, bend right=5] (delta3) to node[below, sloped]{\lab{\rho_{3,4}}} (mult);
    \draw[->, bend left=15] (delta1) to node[below, sloped]{\lab{\sigma_{1,2}}} (m);
    \draw[->, bend left=10] (delta2) to node[below, sloped]{\lab{\sigma_{2,3}}} (m);
    \draw[->, bend left=5] (delta3) to node[below, sloped]{\lab{\sigma_{3,4}}} (m);
    \draw[->, bend right=10] (trr) to node[below, sloped]{\lab{\rho_{1,4}}} (m);
    \draw[->, bend right=10] (mult) to node[below, sloped]{\lab{\rho_{1,4}}} (bll);
  \end{tikzpicture}
  \]
  using the higher product $m_5(\sigma_{3,4},\sigma_{2,3},\sigma_{1,2},x_2,\rho_{1,4})=x_1$.
  See~\cite{LOT2} for more details.

  Note that although the bimodule $\AM(\Id)$ had infinitely many
  nontrivial operations, the bimodule $\DAM(\Id)$ has only finitely
  many. This will be true in general; see the discussion of
  boundedness in~\cite{LOT2}.  (In the terminology there, all of the
  bimodules in this paper are left and right bounded.)
\end{example}

\begin{corollary}\label{cor:enough-to-compute}
  The module $\DAM(\phi)$ is determined by the higher products
  $$m_{m,1,n}(\sigma_{i_1},\dots,\sigma_{i_m},\x,\rho_{j_1},\dots,\rho_{j_n})$$
  on $\AM(\phi)$ where $\sigma_{i_1},\dots,\sigma_{i_m}$ are short
  chords such that $\sigma'_{i_m}\cdots\sigma'_{i_1}\neq 0$.  (In
  particular, the $\sigma_{i_k}$ are distinct and their union is
  connected.)
\end{corollary}
\begin{proof}
  In the bimodule $\DDHalfId\DT\AM(\phi)$, these are the only higher
  products which can lead to non-zero terms in the differential.
\end{proof}

\begin{remark}\label{rmk:DDHalfId-is-dualizer}
  The bimodule $\DDHalfId$ corresponds to the Koszul duality between
  the algebras $\Alg(\PMC)$ and $\Alg(\PMC')$.  See
  \cite[Section~8]{LOTHomPair}.
\end{remark}

\subsection{Equivalence of the two actions}
The reader might wonder if the two actions we have defined are
genuinely different. They are not:
\begin{proposition}\label{prop:actions-equivalent}
  There is an equivalence of categories $\mathcal{F}\co
  \DerBounded(\lsub{\Blg(\PMC)}\ModCat)\to \DerBounded(\lsub{\Clg(\PMC)}\ModCat)$ intertwining
  the actions of the mapping class group of $F(\PMC)$, in the sense
  that the diagram
  \[
  \begin{tikzpicture}
    \node[anchor=mid] at (0,0) (tl) {$\DerBounded(\lsub{\Blg(\PMC)}\ModCat)$};
    \node[anchor=mid] at (5,0) (tr) {$\DerBounded(\lsub{\Blg(\PMC)}\ModCat)$};
    \node[anchor=mid] at (0,-2) (bl) {$\DerBounded(\lsub{\Clg(\PMC)}\ModCat)$};
    \node[anchor=mid] at (5,-2) (br) {$\DerBounded(\lsub{\Clg(\PMC)}\ModCat)$};
    \draw[->] (tl) to node[left] {$\mathcal{F}(\cdot)$} (bl);
    \draw[->] (tr) to node[right] {$\mathcal{F}(\cdot)$} (br);
    \draw[->] (tl) to node[above]{$\cdot\DTP\DAM(\phi)$} (tr);
    \draw[->] (bl) to node[above]{$\cdot\otimes\ADM(\phi)$} (br);
  \end{tikzpicture}
  \]
  commutes.
\end{proposition}
\begin{proof}
  In~\cite{LOT2}, we construct bimodules
  $\lsub{\Blg(\PMC)}\CFDDa(\Id)_{\Clg(\PMC)}$ and 
  $\lsub{\Clg(\PMC)}\CFAAa(\Id)_{\Blg(\PMC)}$ such that for any
  mapping class $\phi$ of $F(\PMC)$,
  \begin{align*}
    \lsub{\Clg(\PMC)}\CFAAa(\Id)_{\Blg(\PMC)}\DTP 
    \lsub{\Blg(\PMC)} \DAM(\phi)_{\Blg(\PMC)}\DTP
    \lsub{\Blg(\PMC)}\CFDDa(\Id)_{\Clg(\PMC)} &\simeq
    \lsub{\Clg(\PMC)}\ADM(\phi)_{\Clg(\PMC)}\\
    \lsub{\Blg(\PMC)}\CFDDa(\Id)_{\Clg(\PMC)}\otimes 
    \lsub{\Clg(\PMC)} \ADM(\phi)_{\Clg(\PMC)}\DTP
    \lsub{\Clg(\PMC)}\CFAAa(\Id)_{\Blg(\PMC)} &\simeq
    \lsub{\Blg(\PMC)}\DAM(\phi)_{\Blg(\PMC)}\\
    \lsub{\Blg(\PMC)}\CFDDa(\Id)_{\Clg(\PMC)}\DTP
    \lsub{\Clg(\PMC)}\CFAAa(\Id)_{\Blg(\PMC)}
    &\simeq \lsub{\Blg(\PMC)}\Blg(\PMC)_{\Blg(\PMC)}\\
    \lsub{\Clg(\PMC)}\CFAAa(\Id)_{\Blg(\PMC)}\DTP
    \lsub{\Blg(\PMC)}\CFDDa(\Id)_{\Clg(\PMC)}
    &\simeq \lsub{\Clg(\PMC)}\Clg(\PMC)_{\Clg(\PMC)}.
  \end{align*}
  So, tensoring with $\CFAAa(\Id_{\PMC})$ gives the desired functor.
\end{proof}

\section{Faithfulness of the action}\label{sec:faithfulness}
To verify that the action is faithful, we start by giving a geometric
interpretation of the rank of~$H_*(M(\phi))$ for $\phi\in\MCG_0(F)$.
By definition, the rank of
$H_*(M(\phi))$ in the idempotent corresponding to $\alpha_i$ and~$\beta_j$
is the Floer homology of $\alpha_i$ with~$\beta_j$.
This Floer homology has a well-known geometric interpretation in terms
of intersection numbers:
\begin{lemma}\label{lem:Floer-is-int-num}
  Let $\alpha$ and $\beta$ be non-isotopic, essential curves in a
  surface $F$, so that $\bdy\alpha\subset \bdy F$,
  $\bdy\beta\subset\bdy F$, $\bdy\alpha\cap\bdy\beta=\emptyset$, and
  $\alpha$ intersects $\beta$ transversely. Let $\HF(\alpha,\beta)$
  denote the Floer homology of the pair $(\alpha,\beta)$. That is,
  $\HF(\alpha,\beta)$ is the homology of the chain complex
  $\CF(\alpha,\beta)$ generated (over $\Field$) by $\alpha\cap\beta$
  and whose differential counts pseudoholomorphic bigons (or,
  equivalently, equivalence classes of immersed bigons) between
  $\alpha$ and $\beta$. Then
  $\dim_\Field(\HF(\alpha,\beta))=i(\alpha,\beta)$.
\end{lemma}

Here $i(\alpha, \beta)$ is geometric intersection number of $\alpha$
and $\beta$: the minimal number of intersections between any two
curves isotopic (relative to the boundary) to $\alpha$ and~$\beta$.
This minimal number is achieved by any curves $\alpha'$ and~$\beta'$
intersecting transversely with no bigons between them.

\begin{proof}
  The Floer homology group $\HF(\alpha,\beta)$ is an isotopy invariant
  of $\alpha$ and $\beta$. If $\alpha'$ and~$\beta'$ are isotopic to
  $\alpha$ and $\beta$ and arranged so that the are no bigons between $\alpha'$
  and $\beta'$ then $\CF(\alpha',\beta')$ has no differential. Thus,
  \[
  \dim_\Field(\HF(\alpha,\beta))=\dim_\Field(\CF(\alpha',\beta'))
  =|\alpha'\cap\beta'|
  =i(\alpha,\beta),
  \]
  as desired.
\end{proof}
To prove faithfulness of the mapping class group action, it suffices to prove:
\begin{theorem}\label{thm:detect-identity}
  The bimodule $\DAM(\phi)$ (respectively $\ADM(\phi)$) is
  quasi-isomorphic to $\lsub{\Blg(\PMC)}\Blg(\PMC)_{\Blg(\PMC)}$
  (respectively $\lsub{\Clg(\PMC)}\Clg(\PMC)_{\Clg(\PMC)}$) if and only
  if $\phi$ is isotopic to the identity.
\end{theorem}
\begin{proof}
  We discuss $\DAM(\phi)$ first.  The functor
  $\Mor_{\Blg(-\PMC')}(\AM(\Id_{\PMC}),\cdot)$ gives an equivalence of
  categories, so it suffices to show that $\AM(\phi)\simeq\AM(\Id)$
  implies $\phi\sim\Id$. Let $I_i$ denote the idempotent corresponding
  to $\eta_i$ (or $\beta_i$) and $J_i$ the idempotent corresponding
  to~$\alpha_i$. Then $I_iH_*(\AM(\phi))J_j$ is the Floer
  homology group $\HF(\beta_i,\alpha_j)$ so, by
  Lemma~\ref{lem:Floer-is-int-num},
  \[
  \dim_\Field I_iH_*(\AM(\phi))J_j=i(\beta_i,\alpha_j).
  \]

  Thus, if $\AM(\phi)\simeq\AM(\Id)$ then $i(\beta_i,
  \alpha_j)=\delta_{i,j}$
  By Lemma~\ref{lem:char-dual-curves}, this implies that
  $\{\beta_i\}$ is isotopic to the set of dual curves $\{\eta_i\}$
  (which are also the $\beta$-curves for the identity map). Thus,
  $\phi$ fixes the curves $\eta_i$ (up to isotopy). Since the
  complement of the $\eta_i$ is a union of disks, and $\phi$ does not
  permute these disks (since $\phi$ fixes the boundary of $F(\PMC)$),
  this implies that $\phi\sim \Id$.

  The statement about $\ADM(\phi)$ follows formally, since
  $\ADM(\phi)\simeq
  \CFAAa(\Id_{\PMC})\DTP\DAM(\phi)\DTP\CFDDa(\Id_{\PMC})$, and tensoring
  with $\CFAAa(\Id)$ and $\CFDDa(\Id)$ give equivalences of categories. Alternatively,
  we can give essentially the same proof as above. Let $\Idem(\PMC)$
  denote the subring of idempotents in $\Clg(\PMC)$. Then
  \[
  I_iH_*(\Idem(-\PMC')\otimes_{\Clg(-\PMC')}\DM(\phi)\otimes_{\Clg(\PMC)}\Idem(\PMC))J_j
\cong \HF(\beta_i,\alpha_j).  
  \]
  (Here, $\Idem(\PMC)$ is a $\Clg(\PMC)$-algebra via the augmentation
  map $\Clg(\PMC)\to\Idem(\PMC)$ sending any non-idempotent element to
  $0$.) The rest of the proof is then the same.
\end{proof}

\begin{proof}[Proof of Theorem~\ref{thm:intro-faithful}]
  This is immediate from Theorem~\ref{thm:detect-identity}, together
  with the identification between $N(\phi)$ and $\CFDAa(\phi,-n/2+1)$
  (Proposition~\ref{prop:bimodules-agree-A}).
\end{proof}

As a corollary, when we iterate a map, the ranks of the homology of
the bimodules grow like the dilatation of a pseudo-Anosov map.

\begin{corollary}
  For $\phi$ a pseudo-Anosov mapping class with dilatation~$\lambda$,
  \[
  \lim_{n \to \infty} \sqrt[n]{\dim_\Field H_*(N(\phi^n))} = 
  \lim_{n \to \infty} \sqrt[n]{\dim_\Field H_*(Q(\phi^n))} = 
  \lambda.
  \]
\end{corollary}

\begin{proof}
  First consider the similar statement for $M(\phi^n)$.  By
  Lemma~\ref{lem:Floer-is-int-num}, 
  \[
  \dim_\Field I_i H_*(M(\phi^n)) J_j = i(\beta_i, \alpha_j) =
  i(\phi^n(\eta_i), \alpha_j).
  \]
  It is a well-known that the intersection numbers in pseudo-Anosov
  maps grow exponentially with the iteration.  More precisely,
  \begin{equation}\label{eq:growth-rate}
  \lim_{n \to \infty} \frac{i(\phi^n(\eta_i), \alpha_j)}{\lambda^n} =
  \mu_s(\eta_i) \mu_u(\alpha_j)
  \end{equation}
  where $\mu_s$ and $\mu_u$ are, respectively, the transverse measures
  on the stable and unstable foliations of $\phi$, suitably
  normalized.  (See, e.g., \cite[Theorem 12.2]{FLP79:TravauxThurston}
  for the theorem for surfaces with no boundary, or \cite[Theorem
  11.5]{FLP79:TravauxThurston} for a related theorem in the case of a
  surface with boundary.)

  For the statement of the corollary, we do not need the precise
  constants on the right-hand side of Equation~\eqref{eq:growth-rate},
  just that they are non-zero.  But $\mu_s(\eta_i) \ne 0$ for any
  pseudo-Anosov map, as otherwise the simple closed curve formed by
  connecting the endpoints of $\eta_i$ along $\bdy F$ would be a
  reducing curve.  (If $\eta_i$ connects two different boundary
  components, consider instead the curve formed by taking two copies
  of $\eta_i$ and connecting the endpoints the long way around~$\bdy
  F$.)  Similarly, $\mu_u(\alpha_j)\neq 0$, so by
  Equation~\eqref{eq:growth-rate}, $i(\phi^n(\eta_i),\alpha_j)$ grows
  as $\lambda^n$. The dimension $\dim_\Field H_*(M(\phi^n))$ is a sum
  of such terms, so $\dim_\Field H_*(M(\phi^n))$ grows as~$\lambda^n$,
  as well.

  By definition, $N(\phi^n) \simeq \Mor_{\Blg(-\PMC')}(M(\Id_{\PMC}),
  M(\phi^n))$.  Since $\Blg(-\PMC')$ and $M(\Id_{\PMC})$ are
  finite\hyp dimensional, $\dim_\Field H_*(N(\phi^n)) \le
  K\dim_\Field H_*(M(\phi^n))$ for some constant~$K$.  Since
  $\Mor_{\Blg(-\PMC')}(M(\Id_{\PMC}), \cdot)$ is an equivalence of
  categories (with inverse given by taking $\Mor$ with another
  bimodule), we also have a similar bound the other direction, proving the
  statement in the corollary for $N(\phi^n)$.

  The statement about $\ADM(\phi)$ follows similarly, since
  $\ADM(\phi)\simeq
  \CFAAa(\Id_{\PMC})\DTP\DAM(\phi)\DTP\CFDDa(\Id_{\PMC})$, and
  both $\CFAAa(\Id_{\PMC})$ and $\CFDDa(\Id_{\PMC})$ are
  finite-dimensional, and tensoring with $\CFAAa(\Id_{\PMC})$
  (respectively $\CFDDa(\Id_{\PMC})$) gives an equivalence of
  categories (where tensoring with $\CFDDa(\Id_{\PMC})$ (respectively
  $\CFAAa(\Id_{\PMC})$) gives the inverse equivalence).
\end{proof}

\begin{remark}
  A similar statement holds if $\phi$ is reducible; then the growth
  rate of the rank of the homology is given by the maximum dilatation
  of any pseudo-Anosov component of $\phi$, as at least one $\alpha_i$
  and one $\eta_j$ must intersect the pseudo-Anosov component.  If
  $\phi$ has no pseudo-Anosov components (i.e., some power of $\phi$ is a composition of
  Dehn twists along pairwise-disjoint curves), the rank of the
  homology grows only linearly.
\end{remark}

\section{Finite generation}\label{sec:finite-generation}
In this section, we briefly review the sense in which the module categories on which the mapping class group is acting are finitely generated.
\begin{definition}
  Given objects $\{M_i\}$ in triangulated category $\Cat$, the
  \emph{subcategory generated by $\{M_i\}$} is the smallest
  triangulated subcategory of $\Cat$ containing all of the~$M_i$. We
  say that \emph{$\{M_i\}$ generate $\Cat$} if the subcategory
  generated by $\{M_i\}$ is, in fact, $\Cat$. We say that \emph{$\Cat$
    is finitely generated} if there is a finite set of objects
  $\{M_i\}$ which generate~$\Cat$.
\end{definition}

Although the definition of finite generation is rather abstract, our
proof that our module categories are finitely generated will be
satisfyingly concrete.
Fix an arc diagram $\PMC$, and let $\Blg=\Blg(\PMC)$, $\Clg=\Clg(\PMC)$.
Before
giving the proof, we develop a little more algebra. Let $\Clg_+$
be the ideal in $\Clg$ generated by all strand diagrams $s$
in which not all strands are horizontal
(so as an $\Field$-vector space, $\Clg$ is the direct sum
of $\Clg_+$ and the subring of idempotents of
$\Clg$). Observe that $\Clg_+$ is nilpotent; for instance,
this follows from the facts that $\Clg$ is finite-dimensional,
the total length gives a grading on~$\Clg$, and $\Clg_+$
is the positively graded part of $\Clg$ with respect to this
grading (compare Remark~\ref{remark:nilpotent}). In particular, for any $\Clg$-module
$M$, the module $\Clg_+\cdot M$ is a proper submodule of~$M$. 

A \emph{simple module} over $\Clg$ is a module $M$ which is
$1$-dimensional over $\Field$ (and so has trivial differential). The
simple modules are in bijective
correspondence with the $2n$ minimal idempotents in~$\Clg$.

\begin{theorem}\label{thm:finite-generation}
  The derived categories $\DerBounded(\lsub{\Blg}\ModCat)$ and $\DerBounded(\lsub{\Clg}\ModCat)$ are finitely generated.
\end{theorem}
\begin{proof}
  We start by proving the statement for
  $\DerBounded(\lsub{\Clg}\ModCat)$; one can give a similar
  proof for $\Blg$, but since we have been working with
  $\Ainf$-modules over $\Blg$ a little extra verbiage is
  required.
  
  We prove that $\Clg$ is generated by the simple modules.
  Our proof is by induction on the dimension over $\Field$ of a
  differential module
  $M\in\DerBounded(\lsub{\Clg}\ModCat)$. There is a short exact
  sequence
  \[
  0\to \Clg_+M\to M\to M/(\Clg_+M)\to 0.
  \]
  Further, $M/\Clg_+M$ is a direct sum of simple modules and $\Clg_+M$
  has strictly smaller dimension than $M$.
  By induction, we can assume
  that $\Clg_+M$ is in the triangulated subcategory generated by the
  simple modules; it follows that $M$ is in this subcategory as well.

  The statement for $\DerBounded(\lsub{\Blg}\ModCat)$ now
  follows from the statement for
  $\DerBounded(\lsub{\Clg}\ModCat)$ and the fact that tensoring
  with $\CFAAa(\Id)$ gives an equivalence between the two categories.
\end{proof}

\begin{remark}
  If we prefer to think of elements of
  $\DerBounded(\lsub{\Clg}\ModCat)$ as projective modules, we
  can give a similar proof using the elementary projective modules
  $\Clg\cdot I$ (for $I$ one of the $n$ minimal idempotents).
\end{remark}

It is not hard to extend the proof of
Theorem~\ref{thm:finite-generation} to give the following:
\begin{theorem}\label{thm:grothendiek}
  The modulo $2$ Grothendieck group $G(\lsub{\Blg}\ModCat)$ of
  differential $\Blg$-modules is isomorphic to
  $H_1(F(\PMC);\ZZ/2)$.  The action of the mapping class group on
  $\lsub{\Blg}\ModCat$ defined in this paper decategorifies to
  the standard action of the mapping class group on
  $H_1(F(\PMC);\ZZ/2)$. The corresponding statements also hold for
  $\lsub{\Clg}\ModCat$, as well as for the Grothendieck groups of
  projective differential modules
  $K_0(\lsub{\Blg}\ModCat)$ and $K_0(\lsub{\Clg}\ModCat)$.
\end{theorem}
\begin{proof}
  The proof of Theorem~\ref{thm:finite-generation} shows that
  the $n$ elementary modules generate
  $G(\lsub{\Clg}\ModCat)$. To see that they are linearly
  independent, consider the algebra map
  \[
  \Clg\to \Clg/\Clg_+=\bigoplus_{i=1}^{n}\Field.
  \]
  This maps the $n$ generators of $G(\lsub{\Clg}\ModCat)$ to
  a basis for
  $G(\lsub{\Clg/\Clg_+}\ModCat)=(\Field)^{n}$.

  To understand the induced mapping class group action, let
  $\alpha_1,\dots,\alpha_{n}$ be the basis of curves for $F(\PMC)$
  specified by the pairs of points in $\Matching$, and let
  $\gamma_1,\dots,\gamma_{n}$ be dual curves. Then, for
  idempotents $I_i$ and $I_j$, the number of generators $g_{ij}$ of $I_i
  \CFDAa(\phi,-n/2+1) I_j$ (as a type \DA\ bimodule) is equal (modulo $2$) to the number of
  intersections between $\phi(\alpha_i)$ and $\gamma_j$. (To see this,
  note that each generator of $\DAM(\phi)$ can be promoted uniquely to
  a generator of $\CFDAa(\phi,-n/2+1)$.) Use the
  $\alpha_i$ to give a basis $[\alpha_1],\dots,[\alpha_{n}]$ for
  $H_1(F(\PMC))$. With respect to this basis,
  \begin{align*}
    \phi_*[\alpha_i]&=(\alpha_i\cdot\gamma_1,\dots,\alpha_i\cdot\gamma_{n})\\
    &\equiv (g_{i,1},\dots,g_{i,n}) \pmod{2}.
  \end{align*}
  This implies that the induced action on Grothendieck groups agrees
  with the action on $H_1(F(\PMC))$.

  The results for $G(\lsub{\Blg}\ModCat)$,
  $K_0(\lsub{\Blg}\ModCat)$ and $K_0(\lsub{\Clg}\ModCat)$
  follow similarly; alternatively, they follow from the fact that all
  of these triangulated categories are equivalent.
\end{proof}

\begin{remark}
  Since we have been working with ungraded differential modules, we
  are forced to use the modulo $2$ Grothendieck groups in
  Theorem~\ref{thm:grothendiek}.
\end{remark}

\begin{remark}
  The proof of Theorem~\ref{thm:grothendiek} immediately extends to
  show that the action of the mapping class group on
  $G(\lsub{\Alg(\PMC)}\ModCat)\cong \Lambda^*H_1(F(\PMC);\ZZ/2)$ is
  the standard action on $\Lambda^*H_1(F(\PMC);\ZZ/2)$.
\end{remark}

\begin{remark}
  In light of Auroux's reformulation of bordered Floer theory in terms
  of partially-wrapped Fukaya categories~\cite{AurouxBordered}, it is
  natural to compare Theorem~\ref{thm:grothendiek} with Abouzaid's
  computation of the Grothendieck group of modules over the Fukaya
  category of a closed surface~\cite{Abouzaid08:FukSurf}: for the
  Fukaya category of a closed surface $F$, the Grothendieck group is
  $H_1(SF;\ZZ)\bigoplus\RR$, where $SF$ is the unit tangent bundle to
  $F$.
\end{remark}
\section{Further questions}\label{sec:further}
The results of this paper suggest several natural questions. Most
prominent among them is whether knowing that the mapping class group
has a faithful representation on a linear category has group-theoretic
consequences. A faithful action of a group on a vector space has many
consequences (like the Tits alternative~\cite{Tits72:alternative} and residual finiteness),
and many of these consequences are known to hold for mapping class groups. It seems plausible
that some of these could be explained by the linear-categorical
actions of the mapping class groups.

A second natural question is whether one can give a similar
linear-categorical action of the mapping class group of a closed
surface.

A question more internal to Heegaard Floer homology is whether the actions on
bordered Floer homology in $\SpinC$-structures between the $(-n/2+1)\st$
and $(n/2-1)\st$ are faithful. It seems likely that they are, but the
techniques of this paper do not apply directly.

Finally, there are many known categorical actions of braid groups. It
would be interesting to know which, if any, of these admit extensions
to mapping class group actions; in particular, this would be a step
towards extending Khovanov-type knot invariants to $3$-manifold
invariants.

\bibliographystyle{hamsalpha}\bibliography{heegaardfloer}

\providecommand{\noopsort}[1]{}
\providecommand{\bysame}{\leavevmode\hbox to3em{\hrulefill}\thinspace}
\providecommand{\href}[2]{#2}
\providecommand{\eprint}{\begingroup \urlstyle{rm}\Url}
\begin{thebibliography}{LOT10b}

\bibitem[Abo08]{Abouzaid08:FukSurf}
Mohammed Abouzaid, \emph{On the {F}ukaya categories of higher genus surfaces},
  Adv. Math. \textbf{217} (2008), no.~3, 1192--1235.

\bibitem[AGW11]{AGW:KSandBordered}
Denis Auroux, J.~Elisenda Grigsby, and Stephan Wehrli, \emph{On
  {K}hovanov-{S}eidel quiver algebras and bordered {F}loer homology}, 2011,
  \eprint{arXiv:1107.2841}.

\bibitem[Aur10]{AurouxBordered}
Denis Auroux, \emph{Fukaya categories of symmetric products and bordered
  {H}eegaard-{F}loer homology}, J. G\"okova Geom. Topol. \textbf{4} (2010),
  1--54, \eprint{arXiv:1001.4323}.

\bibitem[FLP79]{FLP79:TravauxThurston}
A.~Fathi, F.~Laundenbach, and V.~Po{\'e}naru (eds.), \emph{Travaux de
  {T}hurston sur les surfaces}, Ast\'erisque, vol.~66, Soci\'et\'e
  Math\'ematique de France, Paris, 1979, S{\'e}minaire Orsay, with an English
  summary.

\bibitem[GW10]{GrigsbyWehrli:detects}
J.~Elisenda Grigsby and Stephan~M. Wehrli, \emph{On the colored {J}ones
  polynomial, sutured {F}loer homology, and knot {F}loer homology}, Adv. Math.
  \textbf{223} (2010), no.~6, 2114--2165, \eprint{arXiv:0907.4375}.

\bibitem[Kel01]{AinftyAlg}
Bernhard Keller, \emph{Introduction to {$A$}-infinity algebras and modules},
  Homology Homotopy Appl. \textbf{3} (2001), no.~1, 1--35,
  \eprint{arXiv:math.RA/9910179}.

\bibitem[KM11]{KronheimerMrowka11:detect}
Peter~B. Kronheimer and Tomasz Mrowka, \emph{Khovanov homology is an
  unknot-detector}, Publ. Math. Inst. Hautes \'Etudes Sci. (2011), no.~113,
  97--208, \eprint{arXiv:1005.4346}.

\bibitem[KS02]{KhS02:BraidGpAction}
Mikhail Khovanov and Paul Seidel, \emph{Quivers, {F}loer cohomology, and braid
  group actions}, J. Amer. Math. Soc. \textbf{15} (2002), no.~1, 203--271,
  \eprint{arxiv:math.QA/0006056}.

\bibitem[KT07]{KhTomas07:CobordismsCategories}
Mikhail Khovanov and Richard Thomas, \emph{Braid cobordisms, triangulated
  categories, and flag varieties}, Homology Homotopy Appl. \textbf{9} (2007),
  no.~2, 19--94, \eprint{arXiv:math/0609335}.

\bibitem[LOT08]{LOT1}
Robert Lipshitz, Peter~S. Ozsv{\'a}th, and Dylan~P. Thurston, \emph{Bordered
  {H}eegaard {F}loer homology: {I}nvariance and pairing}, 2008,
  \eprint{arXiv:0810.0687v4}.

\bibitem[LOT10a]{LOT2}
\bysame, \emph{Bimodules in bordered {H}eegaard {F}loer homology}, 2010,
  \eprint{arXiv:1003.0598v3}.

\bibitem[LOT10b]{LOT4}
\bysame, \emph{Computing {$\HFa$} by factoring mapping classes}, 2010,
  \eprint{arXiv:1010.2550v3}.

\bibitem[LOT11]{LOTHomPair}
Robert Lipshitz, Peter~S. Ozsv\'ath, and Dylan~P. Thurston, \emph{{H}eegaard
  {F}loer homology as morphism spaces}, Quantum Topology \textbf{2} (2011),
  no.~4, 384--449, \eprint{arXiv:1005.1248}.

\bibitem[Sei02]{Seidel02:FloerMappingClass}
Paul Seidel, \emph{Symplectic {F}loer homology and the mapping class group},
  Pacific J. Math. \textbf{206} (2002), no.~1, 219--229.

\bibitem[Sie11]{Siegel:mcgAction}
Kyler Siegel, \emph{A geometric proof of a faithful linear-categorical surface
  mapping class group action}, 2011, \eprint{arXiv:1108.3676}.

\bibitem[Tit72]{Tits72:alternative}
J.~Tits, \emph{Free subgroups in linear groups}, J. Algebra \textbf{20} (1972),
  250--270.

\bibitem[Zar09]{Zarev09:BorSut}
Rumen Zarev, \emph{Bordered {F}loer homology for sutured manifolds}, 2009,
  \eprint{arXiv:0908.1106}.

\end{thebibliography}
\end{document}